\documentclass[12pt]{article}
\usepackage{amssymb,enumerate,amsthm,amsmath,mathrsfs}
\def\TheDate{May, 2005}\def\TheVersion{1.00}
\def\TheTitle{\ignorespaces
                 Quasi-finite algebras graded by Hamiltonian
\\
                                    and 
\\
                          vertex operator algebras
\\
}\def\TheAuthor{\ignorespaces
Atsushi Matsuo%
\thanks{Graduate School of Mathematical Sciences, University of Tokyo,  Komaba, Tokyo 153-8914, Japan}%
, \ 
Kiyokazu Nagatomo%
\thanks{Department of Pure and Applied Mathematics, Graduate School of Information Science and Technology, Osaka University, Toyonaka, Osaka 560-0043, Japan}%
\ \ and \ 
Akihiro Tsuchiya%
\thanks{Graduate School of Mathematics, Nagoya University,  Furo-cho, Nagoya 464-8602, Japan}
}
\makeatletter
\def\footnoterule{\kern-3\p@\hrule \@width 6cm \kern 2.6\p@}
\def\@makefntext#1{\leftskip1.4em\noindent\llap{${\arabic{footnote}}.\,$}#1}
\def\ps@myheadings{\let\@oddhead\@empty\let\@evenfoot\@empty
\let\@evenhead\@empty\def\@oddfoot{\hfill-- \thepage\ --\hfill}}
\@addtoreset{equation}{section}
\@addtoreset{equation}{subsection}
\newtheoremstyle{theorems}{3ex}{3ex}{\sl}{}{\bf\boldmath}{}{1.0em}{}
\theoremstyle{theorems}
\newtheorem{lemma}{\bfseries Lemma}[subsection]
\newtheorem{proposition}[lemma]{\bfseries Proposition}
\newtheorem{theorem}[lemma]{\bfseries Theorem}
\newtheorem{corollary}[lemma]{\bfseries Corollary}
\newtheoremstyle{definitions}{3ex}{3ex}{\rm}{}{}{}{1.0em}{}
\theoremstyle{definitions}
\newtheorem{definition}[lemma]{\itshape Definition}
\newtheorem{remark}[lemma]{\itshape Remark}
\newtheorem{note}[lemma]{\itshape Note}
\renewcommand\part{
\clearpage
\null\par\vspace{-3.5ex}
\if@noskipsec \leavevmode \fi\par\addvspace{4ex}
\@afterindentfalse\secdef\@part\@spart}
\def\@part[#1]#2{\refstepcounter{part}
{\parindent\z@\raggedright\interlinepenalty\@M\normalfont
\Large \bfseries Part \thepart\par #2
\markright{#1}{}\par}\nobreak\vskip6ex\@afterheading}
\def\@spart#1{}
\makeatother
\begin{document}
\def\markright#1{}
\title{
\vspace{-2cm}{\footnotesize
\hfill
\TheDate
\quad 
{\tt ver.\TheVersion}}\\
\vspace{1cm}%
\TheTitle}
\author{\TheAuthor}\date{}\maketitle
\thispagestyle{empty}
\makeatletter
\def\cases#1{\left\{\,\vcenter{\baselineskip25pt\m@th%
\ialign{$\displaystyle\hfil##\hfil$&\quad##\hfil\crcr#1\crcr}}\right.}
\let\oldbar\bar\def\bar#1{{\bf\oldbar{\fam=-1 #1}}}
\let\oldhat\hat\def\hat#1{{\bf\oldhat{\fam=-1 #1}}}
\def\widebar#1{\!\overline{\!#1\!}\,}\let\tilde\widetilde
\def\eqalign#1
{{\def\normalbaselines{\baselineskip15pt\lineskip4pt
\lineskiplimit4pt}\,\vcenter{\normalbaselines\mathsurround=0pt\ialign{
\hfil$\displaystyle##$\hfil&&$\,$$\displaystyle##$\hfil\crcr\mathstrut
\crcr\noalign{\kern-\baselineskip}#1\crcr\mathstrut\crcr\noalign{\kern
-\baselineskip}}}\,}}
\newcommand{\ensemble}[2]{
\setbox0=\hbox{$#1$#2}\left.\left\{
\vphantom{\vrule height 1.1\ht0 width0pt depth 1.1\dp0}
#1\,\right|\hbox{#2}\right\}}
\let\setdef\ensemble
\makeatother
\newcommand{\ubf}{\mathbf{u}}
\newcommand{\vbf}{\mathbf{v}}
\renewcommand{\emptyset}{\varnothing}
\let\oldcap\cap\renewcommand{\cap}{{\hskip2pt\oldcap\hskip2pt}}
\newcommand{\compo}{
\mathrel{\mathchoice{\scriptstyle\circ}{\scriptstyle\circ}
{\scriptscriptstyle\circ}{\scriptscriptstyle\circ}}}
\def\interior{\mathaccent"7017 }
\newcommand{\NN}{\mathbb{N}}
\newcommand{\ZZ}{\mathbb{Z}}
\newcommand{\CC}{\mathbb{C}}
\newcommand{\KK}{\mathbf{k}}
\newcommand{\Hom}{{\rm Hom}}
\newcommand{\Ker}{{\rm Ker\,}}
\newcommand{\Coker}{{\rm Coker\,}}
\renewcommand{\simeq}{\cong}
\newcommand{\Fmit}{{F}}
\newcommand{\Agraded}{A}
\newcommand{\Agradedhat}{\hat{A}}
\newcommand{\Afinite}{{A}}
\newcommand{\Afiltered}{\mathcal{A}}
\newcommand{\Afilteredv}{\mathcal{A}^\vee}
\def\unity{{1}}
\newcommand{\ham}{\mathit{h}} 
\newcommand{\Hrm}{\mathrm{H}} 
\newcommand{\gap}{{g}}        
\def\Fgraded{\mathrm{F}}
\def\Fgradedv{\mathrm{F}^\vee}
\def\Ffiltered{\mathrm{F}}
\def\Ffilteredv{\mathrm{F}^\vee}
\newcommand\Igraded{\mathrm{I}}
\newcommand\Igradedv{\mathrm{I}^\vee}
\newcommand\Ifiltered{\mathrm{I}}
\newcommand\Ifilteredv{\mathrm{I}^\vee}
\newcommand{\Irm}{\mathrm{I}}
\newcommand{\Jrm}{\mathrm{J}}
\newcommand{\Mmit}{{M}}
\newcommand{\Nmit}{{N}}
\newcommand{\Mcal}{\mathcal{M}}
\newcommand{\Ncal}{\mathcal{N}}
\newcommand{\Qrm}{\mathrm{Q}}
\newcommand{\Qscr}{\mathcal{Q}}
\newcommand{\Pscr}{\mathcal{P}}
\newcommand{\Bregular}{\mathrm{B}}
\newcommand{\Ecal}{\hat{\mathrm{E}}}
\newcommand{\Krm}{\mathrm{K}}
\newcommand{\Prm}{\mathrm{P}}
\newcommand{\Erm}{\mathrm{E}}
\newcommand{\Drm}{\mathrm{D}}
\newcommand{\VOA}{V}          
\newcommand{\vac}{\mathbf{1}} 
\newcommand{\crm}{\mathrm{c}} 
\newcommand{\Crm}{\mathrm{C}} 
\newcommand{\Der}{T} 
\newcommand{\Lvir}{L}
\newcommand{\UEA}{\mathbb{U}} 
\newcommand{\UEAhat}{\hat{\mathbb{U}}} 
\newcommand{\UVOA}{\mathbf{U}} 
\newcommand{\Ucal}{\mathcal{U}}
\newcommand{\Brm}{\mathrm{B}}  
\renewcommand{\Bbb}{\mathbb{B}}
\newcommand{\Grm}{\mathrm{G}}  
\newcommand{\gr}{\mathrm{gr}}  
\newcommand{\Poisson}{{\mathfrak{p}}}
\newcommand{\Sym}{\mathbb{S}}
\newcommand{\SPoisson}{\mathbf{S}}
\newcommand{\SPoissontilde}{\tilde{\SPoisson}}
\newcommand{\agm}{\mathfrak{a}}
\newcommand{\dgm}{\mathfrak{d}}
\newcommand{\ggm}{{\raise1pt\hbox{$\hskip.5pt\mathfrak{g}\hskip.5pt$}}}
\newcommand{\pgm}{\mathfrak{p}}
\newcommand{\qgm}{\mathfrak{q}}
\newcommand{\Dbb}{\mathbb{D}}
%
%
\long\def\neglect#1\tohere{}
%
\neglect
\section{abstract}
\tohere
{\abstract \noindent A general notion of a quasi-finite algebra is introduced as an algebra graded by the set of all integers equipped with topologies on the homogeneous subspaces satisfying certain properties. 
An analogue of the regular bimodule is introduced and various module categories over quasi-finite algebras are described. 
When applied to the current algebras (universal enveloping algebras) of vertex operator algebras satisfying Zhu's $\Crm_2$-finiteness condition, our general consideration derives important consequences on representation theory of such vertex operator algebras.
In particular, the category of modules over such a vertex operator algebra is shown to be equivalent to the category of modules over a finite-dimensional associative algebra. 
}
\section*{Introduction}
%
In order to construct conformal field theories on Riemann surfaces 
associated with a vertex operator algebra $\VOA$ and to obtain their properties such as the finite-dimensionality of the space of conformal 
blocks, factorization of the blocks along the boundaries of the moduli space 
of Riemann surfaces and the fusion functors or the tensor product of
 $\VOA$-modules, we need first to impose an appropriate 
finiteness condition on $\VOA$ and second to study the structure of the 
abelian category of $\VOA$-modules to some extent. 
One of the candidates of such a finiteness condition is the one introduced by 
Y.-C.\ Zhu (\cite{Zhu}), usually called the $\mathrm{C}_2$-finiteness (or
 $\mathrm{C}_2$-cofiniteness), saying that a certain quotient space
 $\VOA/\Crm_2(\VOA)$ is finite-dimensional. 
We will call this condition {\it Zhu's finiteness condition\/} in the rest of the paper.
This condition was used in \cite{Zhu} in an essential way to the proof of 
the modular invariance of characters of $\VOA$-modules, as well as the condition that the category of $\VOA$-modules is semisimple. 
The modular invariance is a part of the characteristic properties of rational conformal field theories. 

Another important ingredient in Zhu's derivation of the modular invariance 
is the use of an associative algebra $A(\VOA)$, called Zhu's algebra, and a 
functor from the category of $\VOA$-modules to the category of
 $A(\VOA)$-modules.
The fundamental result of Zhu is as follows: the functor sends an 
irreducible $\VOA$-module to an irreducible $A(\VOA)$-module in such a way 
that the equivalence classes of irreducibles in the two categories are in 
one-to-one correspondence. In particular, the number of irreducible classes 
is finite if $A(\VOA)$ is finite-dimensional. 
The last property actually follows from Zhu's finiteness condition mentioned above. 

However, the above-mentioned functor need not give us an equivalence of 
categories. 
Specifically, if we include the cases when the category of $\VOA$-modules is 
not semisimple, the algebra $A(\VOA)$ is not enough to understand the 
properties of the category of $\VOA$-modules in general. In this regard, it 
seems to the authors that detailed analysis of the structure of the abelian 
category of $\VOA$-modules has not yet been done. (See \cite{DLM3} for some related results.)

The purpose of the present paper is to fill this gap by utilizing the universal enveloping algebra associated with the vertex operator algebra.

\vskip2ex

The universal enveloping algebra $\UVOA=\UVOA(\VOA)$ associated with $\VOA$ is a certain topological algebra first considered by Frenkel and Zhu in \cite{FZ}.
The presence of topology is inevitable as the universal defining relations of vertex operator algebras contain infinite sums. 
We will call $\UVOA$ the {\it current algebra\/} for short in the rest of the paper. 
More precisely, the current algebra $\UVOA$ is an associative algebra with unity graded by the set of integers equipped with a separated linear topology defined by a certain sequence $\Irm_0,\Irm_1,\ldots$ of open left ideals such that each homogeneous subspace $\UVOA(d)$ is complete with respect to the relative topology. 
Note that the quotient space $\Qrm_n=\UVOA/\Irm_n$ is a discrete space, which inherits a grading from $\UVOA$. We will denote the subspace of degree $d$ by $\Qrm_n(d)$.
The algebra $\UVOA$ is important in that there is a one-to-one correspondence between $\VOA$-modules of certain type and continuous discrete $\UVOA$-modules. (See Section 6 for the precise statement.)
The authors noticed while trying to understand the category 
of $\VOA$-modules by means of $\UVOA$ that the finiteness property should 
better be imposed on $\UVOA$ rather than on $\VOA$ as far as the properties 
of the category as an abelian category are concerned. Thus we arrived at 
the concept of {\it quasi-finiteness}, which is defined as follows: $\UVOA$ 
is quasi-finite if and only if the spaces $\Qrm_n(d)$ are finite-dimensional 
for all $d\in \ZZ$ and $n\in\NN$. 
As we will show in Theorem \ref{136}, Zhu's finiteness condition actually implies the quasi-finiteness of $\UVOA$. 

\vskip2ex

Let us assume that $\UVOA$ is quasi-finite. 
Then the spaces $\Qrm_n(d)$ have generalized eigenspace decompositions with respect to the action of the Virasoro operator $\Lvir_0$, which we will call the Hamiltonian of $\UVOA$.
By using this, we can construct a series of finite-dimensional algebras
 $\UVOA_n$ and functors from the category of continuous discrete
 left $\UVOA$-modules to the category of left $\UVOA_n$-modules which give rise to equivalences of categories when $n$ is sufficiently large. 
Therefore, under the quasi-finiteness, the category of continuous discrete left $\UVOA$-modules as an abelian category is completely described by the properties of the finite-dimensional algebra $\UVOA_n$. 

In application to conformal field theories on Riemann surfaces, we have to enlarge the algebra $\UVOA$ and to take into account the concept of duality of $\UVOA$-modules. 
Therefore we consider two different topologies on $\UVOA$, one is the original left linear one and the other a right linear one, and take the filterwise completion $\Ucal$ and $\Ucal^\vee$ with respect to the two topologies. 
Then we may consider the left $\Ucal$-modules, the right $\Ucal$-modules, the left $\Ucal^\vee$-modules, and the right $\Ucal^\vee$-modules. 
We can now formulate the concept of quasi-finiteness for each of the four and can establish the equivalences or the dualities among them. 
We will call this type of results the {\it finiteness theorems}.

\vskip2ex

As the argument above deriving the equivalences and the dualities of categories only uses quite general properties of $\UVOA$ and $\Lvir_0$, we now postulate them: we will call a topological algebra $\Agraded$ with a distinguished element $\ham$ a {\it quasi-finite algebra graded by Hamiltonain} if it shares the same properties with $\UVOA$ as mentioned above. 
(See Section 2 for the precise definition.)

\vskip2ex

The present paper is divided into two parts: Part I consists in explaining general theory of quasi-finite algebras graded by Hamiltonian and categories of modules over them and Part II in proving that Zhu's finiteness condition on $\VOA$ implies the quasi-finiteness of $\UVOA$. 

\vskip2ex

The plan of Part I is as follows.
We will begin by describing in Section 1 the concept of compatible 
degreewise  topological algebras and the associated topological filtered 
algebras which are modeled on the topological features of $\UVOA$. 
In Section 2, we define the concept of quasi-finite algebras graded by Hamiltonian and give some consequences. 
In particular, we introduce an analogue of the regular bimodule and define certain finite-dimensional algebras $\Afinite_n$. 
We will show that the regular bimodule is a dense subspace of the algebra  (Theorem \ref{073}).
Section 3 is devoted to the equivalence of categories between the continuous 
discrete modules, which we will call {\it exhaustive modules\/}, and the 
category of $\Afinite_n$-modules (Theorem \ref{146}). 
In Section 4, we will introduce a notion of {\it coexhaustive modules\/} 
which is a {\it dual\/} notion of exhaustive modules and establish an 
equivalence between the category of exhaustive modules and the category of 
coexhaustive modules (Theorem \ref{084}). 
The duality of modules will be formulated in Section 5 as the duality between 
quasi-finite exhaustive left modules and quasi-finite coexhaustive right 
modules. 
We will then summarize various equivalences and the dualities for quasi-finite modules (Theorem \ref{233}). 

\vskip2ex

The proof of quasi-finiteness of the current algebra $\UVOA$ in Part II under  Zhu's finiteness condition will be done by considering a certain filtration $\Grm$ on $\UVOA$, which was introduced in \cite{NT}, and a certain universal Poisson algebra $\SPoissontilde=\SPoissontilde(\Poisson)$ associated with $\Poisson=\VOA/\Crm_2(\VOA)$. 
We will construct a surjective homomorphism of Poisson algebra from $\SPoissontilde$ to the degreewise completion $\SPoisson=\SPoisson(\VOA)$ of $\gr^\Grm\UVOA$.
The quasi-finiteness of $\UVOA$ then follows from that of $\SPoissontilde$. 
We will call $\SPoissontilde$ the {\it Poisson current algebra}. 

In Section 6, we will describe the construction and some properties of the current algebra $\UVOA$ associated with a vertex operator algebra $\VOA$.
In Section 7, we will consider the filtration $\Grm$ on $\UVOA$ and show that the associated graded algebra has a structure of a Poisson algebra. 
In Section 8, we will construct the Poisson current algebra
 $\SPoissontilde$ associated with any Poisson algebra $\Poisson$ and we will show that $\SPoissontilde$ is quasi-finite if $\Poisson$ is 
finite-dimensional (Theorem \ref{177}). 
In the final section, we will consider the case when $\Poisson=\VOA/\Crm_2(\VOA)$ and construct the surjective homomorphism $\Psi:\SPoissontilde\rightarrow \SPoisson$ of Poisson algebras (Theorem \ref{176}).
We will then combine these results to show that Zhu's finiteness condition on $\VOA$ implies quasi-finiteness of the current algebra $\UVOA$ (Theorem \ref{136}). 

\vskip2ex

The finiteness theorems mentioned above will give us not only a nice conceptual understanding of the role of Zhu's finiteness condition in representation theory of vertex operator algebras but also a foundation 
in the strategy mentioned at the beginning of this introduction of 
constructing the spaces of conformal blocks on general Riemann surfaces 
and of showing their expected properties. 
This will be developed in our forthcoming paper, which will serve as the continuation of the previous paper \cite{NT} by two of the authors.

\paragraph{\it Acknowledgement.}
The authors thank Dr.\ T.\ Arakawa for discussion. 
A.M.\ thanks Profs.\ T.\ Miwa and S.\ Loktev for conversation. 
K.N.\ thanks Prof.\ N.\ Kawanaka for his continuous encouragements and dedicates the paper to him on the occasion of his sixtieth birthday. 

Results in this paper were presented in part in `Tensor Categories in 
Mathematics and Physics, Vienna, June, 2004', `Workshop on the Moonshine 
Conjectures and Vertex Algebras, Edinburgh, July 2004', `International 
Conference on Infinite Dimensional Lie Theory, Beijing, July 2004', 
`Perspectives arising from vertex algebra theory, Osaka, November 2004'. 
The authors thank the organizers of these conferences. 

This work was partially supported by Grant-in-Aid for Scientific Research 
No.\ 14740005 (Atsushi Matsuo), 16340007 and 16634001 (Kiyokazu Nagatomo) 
and 14204003 (Akihiro Tsuchiya) from the Ministry of Education, Science, 
Sports and Culture, Japan.
%
\part{Quasi-finite algebras graded by Hamiltonian}
%
\section{Linear topologies on graded algebras}
\label{001}
In this section, we introduce compatible degreewise topological algebras and some related notions. 
\subsection{Preliminaries}
\label{197}
We denote the set of integers by $\ZZ$ and the set of nonnegative integers by $\NN$. 

Throughout the paper, we will work over an algebraically closed field $\KK$ of characteristic zero. A vector space always means a vector space over $\KK$ and the scalar multiplication of a vector space is denoted by juxtaposition. 
An algebra means an associative algebra over $\KK$ with unity. 
The multiplication of an algebra $\Agraded$ is denoted by the dot $\cdot$ and the unity by $1_\Agraded$ or simply by $1$. 
For subsets ${S}$ and ${T}$ of $\Agraded$, we denote by ${S}\cdot{T}$ the linear span of the elements of the form $s\cdot t$ with $s\in {S}$ and $t\in {T}$.
The action of $\Agraded$ on an $\Agraded$-module $\Mmit$ is denoted again by the dot. We always assume that the unity $1_\Agraded$ acts by the identity operator. 
We endow the field $\KK$ with the discrete topology. 
Let 
\begin{equation}
  {I}_0, {I}_1,\ldots,{I}_n,\ldots
\end{equation}
be a decreasing sequence of linear subspaces of a vector space ${V}$. 
A {\it linear topology defined by ${I}_n$}\/ means a topology on ${V}$ such that for each $v\in{V}$ the set $\ensemble{v+{I}_n}{$n\in\NN$}$ forms a fundamental system of open neighborhoods of $v$. 
Throughout the paper, a topology on a vector space always means a linear topology given in this way. 
Such a space is usually called a linearly topologized vector space in the literatures.
Let ${V}$ be a vector space with the linear topology defined by ${I}_n$. 
Any subspace $U$ which contains ${I}_n$ for some $n$ is open and closed, and the quotient topology on ${V}/U$ is the discrete topology. 
The closure of a subspace ${U}$ is given by $\bigcap_n\ ({U}+{I}_n)$ and  ${U}$ is dense in ${V}$ if and only if the composite $U\rightarrow V\rightarrow V/I_n$ of canonical maps is surjective for any $n$. 
The completion of ${V}$ with respect to the linear topology defined by ${I}_n$ is the projective limit
\begin{equation}
  \hat{{V}}=\varprojlim_n {V}/{I}_n
\end{equation}
endowed with the projective limit topology, namely the relative topology induced from $\prod_{n=0}^\infty {V}/{I}_n$ with the product topology of the discrete topologies on ${V}/{I}_n$. 
Let $\hat{{I}}_n$ be the closure of the image of ${I}_n$ under the canonical map ${V}\rightarrow \hat{{V}}$.
Then $\hat{{I}}_n$ agrees with the kernel of the canonical map $\hat{{V}}\rightarrow {V}/{I}_n$ and the topology on $\hat{{V}}$ is the linear topology defined by $\hat{{I}}_n$. 
The space ${V}$ is said to be complete if the canonical map ${V}\rightarrow \hat{{V}}$ is a homeomorphism. 
Thus a complete space is separated (i.e.\ Hausdorff) in our convention. 
If ${V}$ is complete then the closure of a subspace $U$ is given by
\begin{equation}
  \hat{U}=\varprojlim_n (U+{I}_n)/{I}_n.
\end{equation}
and for a closed subspace $\Fmit$ the quotient space ${V}/\Fmit$ with the quotient topology is also complete. 
\vskip2ex
We refer the reader to \cite{Bou}, \cite{EGA} and \cite{Mac} for linear topologies. 
\subsection{Compatible degreewise topological algebras}
\label{002}
Let $\Agraded$ be an algebra and suppose given a grading 
\begin{equation}
  \Agraded=\bigoplus_{d\in\ZZ}\Agraded(d)
\end{equation}
indexed by integers such that 
 $
  \Agraded(d)\cdot \Agraded(e)\subset \Agraded(d+e)
$. 
We simply call such an $\Agraded$ a {\it graded algebra\/}.
Let us set
\begin{equation}
  \Fgraded_p\Agraded=\bigoplus_{d\leq p}\Agraded(d)
,\quad
  \Fgraded^\vee_p\Agraded=\bigoplus_{d\geq -p}\Agraded(d)
\end{equation}
where $p$ is an integer.
We call the filtration $\Fgraded$ the {\it associated filtration} and $\Fgraded^\vee$ the {\it opposite filtration}.
Let $\Agraded=\bigoplus_d\Agraded(d)$ be a graded algebra and suppose given a linear topology on each $\Agraded(d)$. 
In such a situation, we will say that $\Agraded$ is endowed with a {\it degreewise topology\/}. 
We assume that the multiplication maps $\Agraded(d)\times\Agraded(e)\rightarrow \Agraded(d+e)$ are continuous. 
We will say that $\Agraded$ is {\it degreewise complete\/} if each $\Agraded(d)$ is complete. 
Now consider the subspace
\begin{equation}
  \Agraded(d)\cap (\Agraded\cdot \Fgraded_{-n-1}\Agraded)
    =\sum_{k\leq -n-1}\Agraded(d-k)\cdot \Agraded(k)
\end{equation}
and let $\Igraded_n(\Agraded(d))$ be its closure in $\Agraded(d)$.
We assume that the sequence $\{\Igraded_n(\Agraded(d))\}$ forms a fundamental system of open neighborhoods of zero in each $\Agraded(d)$.
\begin{definition}
A {\it compatible degreewise topological algebra\/} is a graded algebra $\Agraded$ endowed with a degreewise topology such that all the conditions mentioned above are satisfied. 
A {\it compatible degreewise complete algebra\/} is a compatible degreewise topological algebra $\Agraded$ such that it is degreewise complete. 
\end{definition}
\vskip2ex
Let $\Agraded$ be a compatible degreewise topological algebra. 
Since the multiplication maps $\Agraded(d)\times\Agraded(e)\rightarrow \Agraded(d+e)$ are continuous, we have 
\begin{equation}
\label{225}
    \Agraded(d)\cdot \Igraded_n(\Agraded(e))
      \subset \Igraded_n(\Agraded(d+e))
,\quad
    \Igraded_n(\Agraded(d))\cdot \Agraded(e)
      \subset \Igraded_{n-e}(\Agraded(d+e))
.
\end{equation}
Therefore, for $a\in \Agraded(d)$ and $b\in\Agraded(e)$, we have
\begin{equation}
\label{060}
  (a+\Igraded_{n+e}(\Agraded(d)))\cdot (b+\Igraded_n(\Agraded(e)))
    \subset a\cdot b+\Igraded_n(\Agraded(d+e)).
\end{equation}
%
%

%
\vskip2ex
Instead of $\Igraded_n(\Agraded(d))$, we may consider the closure $\Igraded^\vee_n(\Agraded(d))$ of $\Agraded(d)\cap(\Fgraded^\vee_{-n-1}\Agraded\cdot \Agraded)$.
Since we have
\begin{equation}
\label{207}
  \Igraded^\vee_n(\Agraded(d))=\Igraded_{n-d}(\Agraded(d))
,
\end{equation}
the topology defined by $\Igraded^\vee_n(\Agraded(d))$ agrees with the one defined by $\Igraded_n(\Agraded(d))$.

\vskip2ex

In the sequel, we will use the following terminologies: 
a graded subspace is said to be {\it degreewise dense\/} if each homogenesous subspace is dense; 
the sum of the closures of the homogeneous subspaces of a graded subspace $U$ is called the {\it degreewise closure\/} of $U$. 
\begin{note}
Let $\Agraded$ be a graded algebra endowed with a degreewise topology.
In general, there is no canonical way of extending the topologies on the homogeneous subspaces to the whole space $\Agraded$ which makes $\Agraded$ into a topological algebra. 
\end{note}
\subsection{Degreewise completion}
\label{124}
Let $\Agraded$ be a compatible degreewise topological algebra. Set
\begin{equation}
  \Agradedhat=\bigoplus_d\Agradedhat(d)
\end{equation}
where $\Agradedhat(d)$ is the completion of the space $\Agraded(d)$.
We call $\Agradedhat$ the {\it degreewise completion\/} of $\Agraded$.
Since the multiplication $\Agraded(d)\times\Agraded(e)\rightarrow \Agraded(d+e)$ are continuous, they induce continuous bilinear maps  $\Agradedhat(d)\times\Agradedhat(e)\rightarrow \Agradedhat(d+e)$ which make $\Agradedhat$ into an algebra endowed with a degreewise topology. 
\begin{proposition}
\label{123}
The degreewise completion\/ $\Agradedhat$ is a compatible degreewise complete algebra. 
\end{proposition}
For instance, let $\Agraded=\bigoplus_d\Agraded(d)$ be any graded algebra.
We endow the space $\Agraded(d)$ with the linear topology defined by $\Agraded(d)\cap (\Agraded\cdot\Fgraded_{-n-1}\Agraded)$. 
Then $\Agraded$ becomes a compatible degreewise topological algebra and the degreewise completion $\Agradedhat$ is degreewise complete. 
We call this degreewise topology on $\Agraded$ the {\it standard degreewise topology\/} and the algebra $\Agradedhat$ the {\it standard degreewise completion}.
\begin{note}
Giving topologies on the homogeneous subspaces of a graded algebra in the way described above was considered by some authors, see e.g.\ \cite{FF},  \cite{Ka}, \cite{FZ}, \cite{Mal}, \cite{LW}. 
\end{note}
\subsection{Associated filtered topological algebra}
\label{205}
A topological algebra is said to be left linear if the topology is a linear topology defined by a sequence of left ideals. Right linearity is defined similarly. 
Let $\Agraded$ be a compatible degreewise topological algebra. 
Let us endow the algebra $\Agraded$ with the associated filtration $\Fgraded$.
Consider the left ideal $\Igraded_n(\Agraded)$ of $\Agraded$ defined by
\begin{equation}
  \Igraded_n(\Agraded)=\bigoplus_d\Igraded_n(\Agraded(d)). 
\end{equation}
By the compatibility of $\Agraded$, the space $\Igraded_n(\Agraded)$ is the degreewise closure of $\Agraded\cdot \Fmit_{-n-1}\Agraded$. 
Let us endow the space $\Agraded$ with the topology defined by $\Igraded_n(\Agraded)$. 
Then the relative topology on $\Agraded(d)$ induced from $\Agraded$ agrees with the original topology on $\Agraded(d)$. 
For any $a\in \Fgraded_p\Agraded$ and $b\in\Fgraded_q\Agraded$ we have
\begin{equation}
  (a+\Igraded_{n+q}(\Agraded))\cdot (b+\Igraded_n(\Agraded))\subset a\cdot b+\Igraded_n(\Agraded).
\end{equation}
Hence the multiplication $\Agraded\times\Agraded\rightarrow \Agraded$ is continuous.
In particular, as $\Igraded_n(\Agraded)$ are left ideals, $\Agraded$ is a left linear topological algebra.
We simply call this topology {\it the left linear topology\/} of $\Agraded$.
\vskip2ex
Now consider the opposite filtration $\Fgraded^\vee$ and set
\begin{equation}
  \Igraded^\vee_n(\Agraded)=\bigoplus_d\Igraded^\vee_n(\Agraded(d))
.
\end{equation}
The space $\Igraded^\vee_n(\Agraded)$ is the degreewise closure of $\Fmit^\vee_{-n-1}\Agraded\cdot \Agraded$. 
Then the linear topology on $\Agraded$ defined by $\Igraded^\vee_n(\Agraded)$ now makes $\Agraded$ into a right linear topological algebra. 
We call this topology {\it the right linear topology\/} of $\Agraded$.

\vskip2ex

The left linear topology and the right linear topology on the same algebra $\Agraded$ do {\it not\/} agree with each other in general although the restrictions to $\Agraded(d)$ are the same.
\subsection{Filterwise completion}
\label{003}
Let $\Agraded$ be a compatible degreewise topological algebra. 
Let us endow $\Agraded$ with the left linear topology and let $\Ffiltered_p\Afiltered$ be the completion 
\begin{equation}
  \Ffiltered_p\Afiltered=\varprojlim_n \Fgraded_p\Agraded/\Fgraded_p\Agraded\cap \Igraded_n(\Agraded)
\end{equation}
with the projective limit topology. 
Then the multiplication $\Fgraded_p\Agraded\times \Fgraded_q\Agraded\rightarrow \Fgraded_{p+q}\Agraded$ induces  a multiplication $\Ffiltered_p\Afiltered\times \Ffiltered_q\Afiltered\rightarrow \Ffiltered_{p+q}\Afiltered$.
Consider the space
\begin{equation}
  \Afiltered=\varinjlim_p\Ffiltered_p\Afiltered
\end{equation}
and give it the linear topology defined by 
\begin{equation}
  \Ifiltered_n(\Afiltered)
    =\Ker \bigl(\Afiltered
             \rightarrow \Agraded/\Igraded_n(\Agraded)\bigr)
.
\end{equation}
Then, for any $a\in\Ffiltered_p\Afiltered$ and $b\in\Ffiltered_q\Afiltered$, we have
\begin{equation}
  (a+\Ifiltered_{n+q}(\Afiltered))\cdot (b+\Ifiltered_n(\Afiltered))\subset a\cdot b+\Ifiltered_n(\Afiltered).
\end{equation}
Thus the space $\Afiltered$ becomes a filtered left linear topological algebra such that the subspaces $\Ffiltered_p\Afiltered$ with the relative topologies are complete and that the image of $\Agraded$ under the canonical map $\Agraded\rightarrow \Afiltered$ is a dense subspace of $\Afiltered$.
We call $\Afiltered$ the {\it left linear filterwise completion\/} of $\Agraded$.
\vskip2ex
Now consider the right linear topology on $\Agraded$ and the corresponding  filterwise completion
\begin{equation}
  \Afilteredv
    =\varinjlim_p\Ffilteredv_p\Afilteredv,\quad
     \Ffilteredv_p\Afilteredv
    =\varprojlim_n
     \Fgraded^\vee_p\Agraded/\Fgraded^\vee_p\Agraded\cap\Igraded^\vee_n(\Agraded)
.
\end{equation}
We will call $\Afilteredv$ the {\it right linear filterwise completion\/} of $\Agraded$.
If $\Agraded$ is degreewise complete then the canonical maps
$\Agraded\rightarrow \Afiltered$ and $\Agraded\rightarrow \Afilteredv$ are injective. We will then identify $\Agraded$ with its images.
Thus we have two inclusions
\begin{equation}
\label{221}
  \Afiltered\leftarrow \Agraded\rightarrow \Afilteredv
\end{equation}
such that the relative topologies on each $\Agraded(d)$ induced from $\Afiltered$ and from $\Afilteredv$ agree with the original topology of a compatible degreewise topological algebra. 
Note that the filterwise completions $\Afiltered$ and $\Afilteredv$ are {\it not\/} complete in general. 
\subsection{Hamiltonian of a graded algebra}
Let $\Agraded$ be a graded algebra.
An element $\ham\in \Agraded$ is called a {\it Hamiltonian\/} of $\Agraded$ if
\begin{equation}
  \Agraded(d)=\ensemble{a\in\Agraded}{$[\ham,a]=da$}
\end{equation}
holds for any $d$, where $[\ham,a]=\ham\cdot a-a\cdot \ham$ denotes the commutator.
If this is the case then any central element as well as $\ham$ itself belongs to $\Agraded(0)$ and an element $\ham'$ is a Hamiltonian if and only if $\ham-\ham'$ is central.
An {\it algebra graded by Hamiltonian\/} is a pair $(\Agraded,\ham)$ of a graded algebra $\Agraded$ and a Hamiltonian $\ham\in \Agraded(0)$. 
We will denote by $\Hrm$ the subalgebra of $\Agraded(0)$ generated by the Hamiltonian $\ham$. 
\begin{remark}
Let $(\Agraded,\ham)$ be a compatible degreewise complete algebra graded by Hamiltonian. 
Then the images of the canonical injections (\ref{221}) are given respectively by
 $
  \sum_d\Afiltered(d)
$ 
and
 $
  \sum_d\Afilteredv(d)
$, 
where
\begin{equation}
  \Afiltered(d)=\ensemble{a\in\Afiltered}{$[\ham,a]=da$}
,\quad
  \Afilteredv(d)=\ensemble{a\in\Afilteredv}{$[\ham,a]=da$}.
\end{equation}
\end{remark}
In the rest of Part I, we will be mainly concerned with a compatible degreewise complete algebra with Hamiltonian. 
%
%
\section{Quasi-finite algebras}
\label{122}
In this section, we will formulate a finiteness condition on compatible degreewise complete algebras and describe its consequences. 
In particular, we will formulate an analogue of the regular bimodule as a degreewise dense subspace of the algebra. 
\subsection{Canonical quotient modules}
\label{218}
Let $\Agraded$ be a compatible degreewise topological algebra and recall the spaces $\Igraded_n(\Agraded)$. 
We set
\begin{equation}
  \Qrm_n=\Agraded/\Igraded_n(\Agraded)
.
\end{equation}
Since $\Igraded_n(\Agraded)$ is a left ideal, the quotient $\Qrm_n$ is a left $\Agraded$-module. 
We will call this module the {\it left canonical quotient module}.
By (\ref{225}), we have $\Igraded_n(\Agraded)\cdot \Fgraded_0\Agraded\subset \Igraded_n(\Agraded)$.
Hence the multiplication $\Agraded\times \Fgraded_0\Agraded\rightarrow \Agraded$ induces an action of $\Fgraded_0\Agraded$ on $\Qrm_n$ for any $n$. 
Thus the space $\Qrm_n$ is an $(\Agraded,\Fgraded_0\Agraded)$-bimodule. 
In particular, it is an $(\Agraded(0),\Agraded(0))$-bimodule. 
Since $\Igraded_n(\Agraded)$ is a graded subspace, the grading of $\Agraded$ induces a grading on the quotient by setting
 $
  \Qrm_n(d)=\Agraded(d)/\Igraded_n(\Agraded(d))
$.
If $d\leq -n-1$ then since $\Agraded(d)\subset \Igraded_n(\Agraded)$ we have  $\Qrm_n(d)=0$. 
Thus
\begin{equation}
  \Qrm_n=\bigoplus_{d=-n}^\infty \Qrm_n(d)
.
\end{equation}
\begin{lemma}
\label{226}
For any $v\in \Qrm_n$ there exists an $m$ such that\/ $\Igraded_m(\Agraded)\cdot v=0$. 
\begin{proof}
By (\ref{225}) we have $\Igraded_{n+d}(\Agraded)\cdot \Agraded(d)\subset \Igraded_n(\Agraded)$ and hence $
\Igraded_{m}(\Agraded)\cdot \Qrm_n(d)=0
$ for $m=n+d$. 
\end{proof}
\end{lemma}
Note that the action $\Agraded\times \Qrm_n\rightarrow \Qrm_n$ is continuous when $\Agraded$ is endowed with the left linear topology and $\Qrm_n$ with the discrete topology. 

\begin{lemma}
\label{227}
Let $v$ be an element of\/ $\Qrm_n$.
Then $\Fgraded_{-n-1}\Agraded\cdot v=0$ implies\/ $\Igraded_n(\Agraded)\cdot v=0$. 
\begin{proof}
Suppose $\Fgraded_{-n-1}\Agraded\cdot v=0$ and choose an $m$ such that $\Igraded_m(\Agraded)\cdot v=0$. 
Since $\Igraded_n(\Agraded)$ is the closure of $\Agraded\cdot \Fgraded_{-n-1}\Agraded$, we have $\Igraded_n(\Agraded)\subset \Agraded\cdot \Fgraded_{-n-1}\Agraded+\Igraded_m(\Agraded)$. 
Hence $\Igraded_n(\Agraded)\cdot v\subset \Agraded\cdot \Fgraded_{-n-1}\Agraded\cdot v+\Igraded_m(\Agraded)\cdot v=0$.
\end{proof}
\end{lemma}
We may likewise consider the right canonical quotient module $\Qrm^\vee_n=\Agraded/\Igraded^\vee_n(\Agraded)$. 
We have analogous statements for the right canonical quotient modules as well. 
By the definitions of the filterwise completions $\Afiltered$ and
 $\Afilteredv$, the spaces $\Qrm_n$ and $\Qrm^\vee_n$ are canonically 
isomorphic to the spaces $\Afiltered/\Igraded_n(\Afiltered)$ and
 $\Afilteredv/\Igraded^\vee_n(\Afilteredv)$, respectively.
\subsection{Quasi-finite algebra graded by Hamiltonian}
\label{014}
Let $\Agraded$ be a compatible degreewise topological algebra. 
Let us consider the subspace $\Agraded(0)$ of degree zero, which is a subalgebra of $\Agraded$. 
Then $\Igraded_n(\Agraded(0))$ and $\Igraded^\vee_n(\Agraded(0))$ agree with each other and they give a two-sided ideal of $\Agraded(0)$. 
Therefore, the quotient $\Qrm_n(0)=\Agraded(0)/\Igraded_n(\Agraded(0))$ is an algebra.

\begin{lemma}
\label{235}
The space\/ $\Qrm_n(d)$ is a\/ $(\Qrm_{n+d}(0),\Qrm_n(0))$-bimodule. 
\begin{proof}
By (\ref{225}), we have 
$\Igraded_{n+d}(\Agraded(0))\cdot \Agraded(d)\subset \Igraded_n(\Agraded(d))$ and $\Agraded(d)\cdot \Igraded_n(\Agraded(0))\subset \Igraded_n(\Agraded(d))$. 
The result follows. 
\end{proof}
\end{lemma}

Let $\ham$ be a Hamiltonian and let $\ham_{n}$ be the image of $\ham$ in the quotient $\Qrm_n(0)$.
Let $\Hrm_n$ be the subalgebra of $\Qrm_n(0)$ generated by $\ham_{n}$. 
By Lemma \ref{235}, the space $\Qrm_n(d)$ is in particular an $(\Hrm_{n+d},\Hrm_n)$-bimodule. 

\begin{definition}
\label{222}
A {\it weakly quasi-finite algebra\/} is a compatible degreewise complete algebra such that $\Qrm_n(0)$ are finite-dimensional for all $n$. 
A {\it weakly quasi-finite algebra graded by Hamiltonian\/} is a pair $(\Agraded,\ham)$ of a weakly quasi-finite algebra $\Agraded$ and a Hamiltonian $\ham$ of $\Agraded$.
\end{definition}
\noindent

Let $(\Agraded,\ham)$ be a weakly quasi-finite algebra graded by Hamiltonian and let $\ham_{n}$ and $\Hrm_n$ be as above. 
Then $\Hrm_n$ is a finite-dimensional commutative algebra for any $n\in\NN$.
This last property is what we will need in practice in considering a weakly quasi-finite algebra graded by Hamiltonian. 

\vskip2ex
Now let us consider the following conditions for each $d$:
\vskip2ex
\begin{enumerate}[(a)]
\item 
The spaces $\Qrm_n(d)$ are finite-dimensional for all $n$.
\vskip1ex
\item
The spaces $\Qrm^\vee_n(d)$ are finite-dimensional for all $n$.
\end{enumerate}
\vskip1ex
\noindent
Thanks to the relation (\ref{207}), these conditions are actually equivalent.
\begin{definition}
\label{090}
A {\it quasi-finite algebra\/} is a compatible degreewise complete algebra  such that the equivalent conditions (a) and (b) above are satisfied for all integers $d$.
A {\it quasi-finite algebra graded by Hamiltonian\/} is a pair $(\Agraded,\ham)$ of a quasi-finite algebra $\Agraded$ and a Hamiltonian $\ham$ of $\Agraded$.
\end{definition}
\subsection{Spectrum of the Hamiltonian}
\label{193}
Let $(\Agraded,\ham)$ be a weakly quasi-finite algebra graded by Hamiltonian and recall the image $\ham_{n}$ of $\ham$ in $\Qrm_n(0)$.
We let $\varphi_n$ be the minimal polynomial of $\ham_{n}$ and let $\Omega_n$ be the set of roots of $\varphi_n$. 
Then the set $\Omega_n$ agrees with the eigenvalues of the left action of $\ham$ on $\Qrm_n(0)$. 
Let us introduce a partial order on $\KK$ by letting $\lambda\succcurlyeq \mu$ when $\lambda-\mu$ is a nonnegative integer and let $\Gamma_0$ be the set of minimal elements of $\Omega_0$. 
We set 
\begin{equation}
  \Gamma_\infty
  =\ensemble{\gamma+k}{$\gamma\in\Gamma_0$ and $k\in\NN$}
.
\end{equation}
We will also use the following notation: 
\begin{equation}
  \Gamma_m
  =\ensemble{\gamma+k}{$\gamma\in\Gamma_0$ and $k\in\NN$ with $0\leq k\leq m$}
.
\end{equation}
Let $\gap$ be the maximum of the integral differences among elements of $\Omega_0$:
\begin{equation}
\label{231}
  \gap=\max\ensemble{\lambda-\mu}{$\lambda,\mu\in\Omega_0$, $\lambda-\mu\in\NN$}
.
\end{equation}
Then we have $\Gamma_0\subset \Omega_0\subset \Gamma_\gap$.

\vskip2ex

Recall the subalgebras $\Hrm\subset \Agraded(0)$ and $\Hrm_n\subset \Qrm_n(0)$.
For any left $\Hrm$-module ${W}$, we denote by ${W}[\lambda]$ the generalized eigenspace of the Hamiltonian $\ham$ acting on ${W}$:
\begin{equation}
  {W}[\lambda]=\ensemble{v\in{W}}{$(h-\lambda)^r\cdot v=0$ for some $r\in\NN$}.
\end{equation}
If ${W}$ is an $\Hrm_n$-module then ${W}[\lambda]\ne 0$ implies $\lambda\in \Omega_n$. 

\vskip2ex

Now consider the space 
\begin{equation}
  \Krm_m(\Qrm_n)
  =\ensemble{v\in\Qrm_n}{$\Igraded_m(\Agraded)\cdot v=0$}.
\end{equation}
Then this is a left $\Qrm_m(0)$-module and hence a left $\Hrm_m$-module. 
Thus the set of the eigenvalues of the left action of $\ham$ on $\Krm_m(\Qrm_n)$ is contained in the set $\Omega_m$. 
Note that $v\in \Krm_m(\Qrm_n)$ if and only if $\Fgraded_{-m-1}\Agraded\cdot v=0$ by Lemma \ref{227}.
\begin{proposition}
\label{036}
The set\/ $\Omega_n$ is a subset of\/ $\Gamma_{n+\gap}$.
\begin{proof}
Let $\lambda$ be an element of $\Omega_n$.
Then there exists a generalized eigenvector $v$ in $\Qrm_{n}(0)\subset \Krm_n(\Qrm_n)$ with the eigenvalue $\lambda$.
Since $v\ne 0$ and $\Fgraded_{-n-1}\Agraded\cdot v=0$, there exists a nonnegative integer $k$ with $0\leq k\leq n$ for which $v\notin \Krm_{k-1}(\Qrm_n)$ but $v\in\Krm_k(\Qrm_n)$. 
Then $\Agraded(-k)\cdot v$ is a nonzero subspace of $\Krm_0(\Qrm_n)$.
Since $\Agraded(-k)\cdot v\subset \Qrm_n[\lambda-k]$, we have $\lambda-k\in \Omega_0$ and hence $\lambda\in \Gamma_{k+\gap}\subset \Gamma_{n+\gap}$.
\end{proof}
\end{proposition}
Now Lemma \ref{226} implies $\Qrm_n=\bigcup_{m=0}^\infty \Krm_m(\Qrm_n)$. 
Therefore, Proposition \ref{036} implies that 
\begin{equation}
  \Qrm_n=\bigoplus_{\lambda\in\Gamma_\infty}\Qrm_n[\lambda].
\end{equation}
\begin{lemma}
\label{228}
The space\/ $\Qrm_n[\lambda]$ with $\lambda\in\Gamma_m$ is a subspace of $\Krm_{m}(\Qrm_n)$.
\begin{proof}
Since $\lambda-m-1\notin \Gamma_\infty$, we have $\Fgraded_{-m-1}\Agraded\cdot \Qrm_n[\lambda]=0$. 
Hence $\Qrm_n[\lambda]\subset \Krm_m(\Qrm_n)$. 
\end{proof}
\end{lemma}
Let $\ell$ be the maximum of the multiplicities of the roots of the minimal polynomial $\varphi_\gap$. 
\begin{proposition}
\label{035}
The multiplicities of the roots of $\varphi_n$ are at most $\ell$ for any nonnegative integer $n$.
\begin{proof}
Let ${v}$ be an element of $\Qrm_n$. 
We show that $(h-\lambda)^k\cdot v=0$ for some $k$ implies $(h-\lambda)^\ell\cdot v=0$.
Let $\lambda$ be a minimal counterexample to this claim: there exists a nonzero vector $v$ such that $(h-\lambda)^k\cdot v=0$ for some $k$ but $(h-\lambda)^{\ell}\cdot v\ne 0$.
Then, by the minimality, we have $\Fgraded_{-1}\Agraded\cdot (\ham_{n}-\lambda)^{\ell}\cdot v=0$ and hence $(\ham_{n}-\lambda)^{\ell}\cdot v\in \Krm_0(\Qrm_n)$ by Lemma \ref{227}. 
Hence $\lambda\in\Omega_0\subset \Gamma_\gap$ and so $v\in\Qrm_n[\lambda]\subset \Krm_\gap(\Qrm_n)$ by Lemma \ref{228}.
Hence $(h-\lambda)^\ell\cdot v=0$, a contradiction.
\end{proof}
\end{proposition}
\subsection{The regular bimodule}
\label{065}
Let $(\Agraded,\ham)$ be a weakly quasi-finite algebra graded by Hamiltonian. 
We denote by $\Agraded[\lambda,\mu]$ the simultaneous generalized eigenspace of the left and the right actions of the Hamiltonian $\ham$ on $\Agraded$:
\begin{equation}
  \Agraded[\lambda,\mu]
    =\ensemble{a\in\Agraded}%
              {$(h-\lambda)^r\cdot a=a\cdot (h-\mu)^r=0$ for some $r$}.
\end{equation}
We will call this $\lambda$ the {\it left eigenvalue\/} and $\mu$ the {\it right eigenvalue\/} of $\ham$. 

By the definition, $\Agraded[\lambda,\mu]\ne 0$ implies that $d=\lambda-\mu$ is an integer and $\Agraded[\lambda,\mu]\subset\Agraded(d)$. 
We also note that
\begin{equation}
\label{219}
  \Agraded(d)\cdot \Agraded[\lambda,\mu]\subset \Agraded[\lambda+d,\mu],\quad
  \Agraded[\lambda,\mu]\cdot \Agraded(e)\subset \Agraded[\lambda,\mu-e].
\end{equation}
It is easy to see that
\begin{equation}
\hbox{
 $
  \Agraded[\kappa,\lambda]\cdot \Agraded[\mu,\nu]=0
 $
if $\lambda\ne \mu$.
}
\end{equation}
Indeed, if $a\cdot (h-\lambda)^k=0$ and $(h-\mu)^m\cdot b=0$ with $\lambda\ne \mu$ then $a\cdot b=a\cdot 1\cdot b=a\cdot (h-\lambda)^kf(h)\cdot b+a\cdot (h-\mu)^mg(h)\cdot b=0$ for some polynomials $f(x)$ and $g(x)$.
Also note that 
\begin{equation}
\label{208}
  \Agraded[\kappa,\lambda]\cdot \Agraded[\lambda,\mu]\subset \Agraded[\kappa,\mu]
.
\end{equation}
We now set
\begin{equation}
  \Bregular=\sum_d\Bregular(d),\quad 
  \Bregular(d)=\sum_{\lambda-\mu=d}\Agraded[\lambda,\mu].
\end{equation}
Then by (\ref{219}) this is an $(\Agraded,\Agraded)$-bimodule. 
We call $\Bregular$ the {\it regular bimodule\/} of $\Agraded$.
\begin{remark}
\label{182}
The structure of an $(\Agraded,\Agraded)$-bimodule on $\Bregular$ actually prolongs to a structure of an $(\Afiltered,\Afilteredv)$-bimodule.
(See Proposition \ref{098}.)
\end{remark}
\subsection{Denseness of the regular bimodule}
\label{212}
To investigate the spaces $\Agraded[\lambda,\mu]$ we consider the left canonical quotient module $\Qrm_n$, which is an $(\Agraded,\Fgraded_0\Agraded)$-bimodule.
Consider the simultaneous generalized eigenspaces of the left and the right actions of $\ham$ on $\Qrm_n$:
\begin{equation}
  \Qrm_n[\lambda,\mu]
    =\ensemble{v\in\Qrm_n}%
              {$(h-\lambda)^r\cdot v=v\cdot (h-\mu)^r=0$ for some $r$}.
\end{equation}
Then $\Qrm_n[\lambda,\mu]\ne 0$ implies that $d=\lambda-\mu$ is an integer and that $\Qrm_n[\lambda,\mu]\subset\Qrm_n(d)$. 

Recall that the space $\Qrm_n(d)$ is a $(\Qrm_{n+d}(0),\Qrm_n(0))$-bimodule. 
In particular, it is an $(\Hrm_{n+d},\Hrm_n)$-bimodule. 
Hence we have
\begin{equation}
  \Qrm_n(d)
   =\sum_{\mu\in\Omega_n}\Qrm_n[\mu+d,\mu].
\end{equation}
Note that $\Agraded(d)\cdot \Qrm_n[\lambda,\mu]\subset \Qrm_n[\lambda+d,\mu]$ and that $\Qrm_n[\lambda,\mu]\cdot \Agraded(e)\subset \Qrm_n[\lambda,\mu-e]$ whenever $e\leq 0$.
\begin{lemma}
\label{075}
If\/ $\Agraded[\lambda,\mu]\ne 0$ then $\lambda,\mu\in\Gamma_\infty$.
\begin{proof}
Let $a$ be a nonzero element of $\Agraded[\lambda,\mu]$ with $d=\lambda-\mu$. 
Since we have assumed that $\Agraded(d)$ is complete, it is separated.
Hence there exists an $n$ such that the image $v$ of $a$ in $\Qrm_n(d)$ is nonzero. 
Then since $v\in\Qrm_n[\lambda,\mu]$, we have $\lambda\in\Gamma_\infty$ and $\mu\in\Gamma_n\subset \Gamma_\infty$.
\end{proof}
\end{lemma}
\begin{lemma}
\label{074}
If $\mu\in \Gamma_m$ and $n\geq m$ then the restriction\/ 
 $ 
  \Qrm_n[\lambda,\mu]
    \rightarrow \Qrm_m[\lambda,\mu]
 $ 
of the canonical surjection is an isomorphism for any $\lambda$.
\begin{proof}
Set $d=\lambda-\mu$.
Consider the kernel $\Igraded_m(\Agraded)/\Igraded_n(\Agraded)$ of the canonical surjection $\Qrm_n\rightarrow \Qrm_m$.
Consider the map
 $
  \Agraded\times\Fgraded_{-m-1}\Agraded\rightarrow \Agraded
 $
which induces
\begin{equation}
  \Agraded\times \Fgraded_{-m-1}\Agraded
        \rightarrow \Igraded_{m}(\Agraded)/\Igraded_n(\Agraded).
\end{equation}
Since $\Igraded_m(\Agraded)$ is the closure of $\Agraded\cdot \Fgraded_{-m-1}\Agraded$, this map is surjective.
Hence the induced map
\begin{equation}
  \Qrm_n\times \Fgraded_{-m-1}\Agraded
        \rightarrow \Igraded_{m}(\Agraded)/\Igraded_n(\Agraded)
\end{equation}
is also surjective.
The right eigenvalues of $\ham$ on $\Qrm_n\cdot \Fgraded_{-m-1}\Agraded$ exceed those of $\Qrm_n$ by more than $m$. 
Hence it follows that the set of the right eigenvalues on the space $\Igraded_{m}(\Agraded)/\Igraded_n(\Agraded)$ does not intersect $\Gamma_m$. 
This implies the result.
\end{proof}
\end{lemma}

By this lemma, we have the following result.

\begin{proposition}
\label{230}
Let\/ $(\Agraded,\ham)$ be a weakly quasi-finite algebra graded by Hamiltonian and let $m$  be a nonnegative integer.
Then the canonical surjection\/ $\Agraded[\lambda,\mu]\rightarrow \Qrm_n[\lambda,\mu]$ is an isomorphism whenever $\mu\in\Gamma_m$ and $n\geq m$. 
\begin{proof}
By Lemma \ref{074}, we have a canonical splitting 
 $
  \Qrm_m[\lambda,\mu]\rightarrow \varprojlim_n\Qrm_n[\lambda,\mu]
$. 
Set $d=\lambda-\mu$.
Since $\Agraded(d)$ is complete, the projective limit is isomorphic to the subspace $\Agraded[\lambda,\mu]$ of $\Agraded(d)$.
\end{proof}
\end{proposition}

The following is the first main result of Part I. 
\begin{theorem}
\label{073}
Let\/ $(\Agraded,\ham)$ be a weakly quasi-finite algebra graded by Hamiltonian.
Then the regular bimodule\/ $\Bregular$ is degreewise dense in\/ $\Agraded$. 
\begin{proof}
By Proposition \ref{230}, the map 
$
  \Bregular(d)\rightarrow \Qrm_n(d)=\Agraded(d)/\Igraded_n(\Agraded(d))
$
is surjective for any $n$. 
This means that $\Bregular(d)$ is dense in $\Agraded(d)$. 
\end{proof}
\end{theorem}

Note that the whole space $\Bregular$ is dense in $\Agraded$ with respect to the left linear and the right linear topologies of $\Agraded$.  
\subsection{The associated finite algebras}
\label{100}
Let $(\Agraded,\ham)$ be a weakly quasi-finite algebra graded by Hamiltonian. 
We set
\begin{equation}
  \Afinite_n
    =\sum_{\lambda,\mu\in\Gamma_n}\Agraded[\lambda,\mu].
\end{equation}
Then it follows from (\ref{208}) that the space $\Afinite_n$ is closed under the multiplication of $\Agraded$. 
By Proposition \ref{230}, we may identify $\Afinite_n(0)$ with a subspace of $\Qrm_n(0)$.
Recall the elements $1_{n}$ and $\ham_{n}$ of $\Qrm_n(0)$, which are the images of $1$ and $\ham$ under the map $\Agraded(0)\rightarrow \Qrm_n(0)$, respectively.
Let 
\begin{equation}
  1_{n}=\sum_{\lambda\in\Omega_n}1_{n}[\lambda,\lambda]
,\quad 
  \ham_{n}=\sum_{\lambda\in\Omega_n}\ham_{n}[\lambda,\lambda]
\end{equation}
be the decompositions to sums of simultaneous generalized eigenvectors.
We set
\begin{equation}
  {1}_{\Afinite_n}=\sum_{\lambda\in\Gamma_n}1_{n}[\lambda,\lambda]
,\quad 
  {\ham}_{\Afinite_n}=\sum_{\lambda\in\Gamma_n}\ham_{n}[\lambda,\lambda],
\end{equation}
and regard them as elements of $\Afinite_n$.
\begin{proposition}
\label{101}
The graded algebra structure on\/ $\Agraded$ induces a graded algebra structure on\/ $\Afinite_n$ for which\/ ${1}_{\Afinite_n}$ is the unity and\/ ${\ham}_{\Afinite_n}$ is a Hamiltonian.
\end{proposition}
\vskip2ex
%

%
%
%
%
%
%
We set
\begin{equation}
  \Prm_{n}
    =\sum_{\lambda\in\Gamma_\infty}\sum_{\mu\in\Gamma_n}\Agraded[\lambda,\mu]
,\quad
  \Prm^\vee_{m}
    =\sum_{\lambda\in\Gamma_m}\sum_{\mu\in\Gamma_\infty}\Agraded[\lambda,\mu]
.
\end{equation}
Then $\Prm_{n}$ is an $(\Agraded,\Afinite_n)$-bimodule and  $\Prm^\vee_{m}$ is an $(\Afinite_m,\Agraded)$-bimodule, which  will play prominent roles in the next section.
\begin{remark}
\label{220}
The structure of an $(\Agraded,\Afinite_n)$-bimodule on $\Prm_n$ prolongs to a structure of an $(\Afiltered,\Afinite_n)$-bimodule and the structure of an $(\Afinite_m,\Agraded)$-bimodule on $\Prm^\vee_m$ to an $(\Afinite_m,\Afilteredv)$-bimodule. (See Remark \ref{182} and Proposition \ref{098}.)
\end{remark}
Let us consider the case when $\Agraded$ is quasi-finite. 
\begin{proposition}
\label{232}
Let\/ $(\Agraded,\ham)$ be a quasi-finite algebra graded by Hamiltonian.
Then\/ $\Afinite_n$ are finite-dimensional for all $n$.
\begin{proof}
For each $\lambda$ and each $\mu$, the space $\Agraded[\lambda,\mu]$ is isomorphic to $\Qrm_m[\lambda,\mu]$ for sufficiently large $m$ by Proposition \ref{230}. 
Since the range $\Gamma_n$ of $\lambda$ and $\mu$ is finite and $\Qrm_m[\lambda,\mu]\subset \Qrm_m(\lambda-\mu)$ is finite-dimensional, the space $\Afinite_n$ is also finite-dimensional.
\end{proof}
\end{proposition}
%
%
\section{Exhaustive modules}
\label{008}
We will define the notion of exhaustive modules and investigate the properties of the category of such modules over a weakly quasi-finite algebra $\Agraded$ with Hamiltonian. 
In particular, we will show that the category of exhaustive $\Agraded$-modules and the category of $\Afinite_n$-modules are equivalent if $n\geq \gap$. 
\subsection{Exhaustive modules}
\label{004}
Let $(\Agraded,\ham)$ be a weakly quasi-finite algebra graded by Hamiltonian and endow it with the left linear topology. 
\begin{definition}
A left $\Agraded$-module $\Mmit$ is {\it exhaustive\/} if for any $v\in \Mmit$ there exists an $m$ such that $\Igraded_m(\Agraded)\cdot v=0$. 
\end{definition}
By this definition, it is easy to see that a left $\Agraded$-module $\Mmit$ is exhaustive if and only if the action $\Agraded\times \Mmit\rightarrow \Mmit$ is continuous with respect to the left linear topology on $\Agraded$ and the discrete topology on $\Mmit$.
Note, however, that a topology on an exhaustive left $\Agraded$-module $\Mmit$  need {\it not\/} be the discrete topology in order for the action $\Agraded\times \Mmit\rightarrow \Mmit$ to be continuous with respect to the left linear topology on $\Agraded$. 
\vskip2ex
\begin{lemma}
\label{106}
Let\/ $\Mmit$ be an exhaustive left\/ $\Agraded$-module and let $v$ be an element of\/ $\Mmit$. 
Then\/ $\Fgraded_{-n-1}\Agraded\cdot v=0$ implies\/ $\Igraded_{n}(\Agraded)\cdot v=0$.
\begin{proof}
See the proof of Lemma \ref{227}. 
\end{proof}
\end{lemma}
Let us now consider the left linear filterwise completion $\Afiltered$. 
We may analogously define the notion of exhaustive left $\Afiltered$-modules as follows.
\begin{definition}
A left $\Afiltered$-module $\Mmit$ is {\it exhaustive\/} if for any $v\in \Mmit$ there exists an $m$ such that $\Igraded_m(\Afiltered)\cdot v=0$. 
\end{definition}
It is fairly clear that results similar to those given above hold for exhaustive left $\Afiltered$-modules. 
The following proposition is a particular case of a general fact on topological algebras. (See e.g.\ \cite{Bou}.)
\begin{proposition}
\label{098}
For any exhaustive left\/ $\Agraded$-module $\Mmit$, the action\/ $\Agraded\times \Mmit\rightarrow \Mmit$ induces an exhaustive left $\Afiltered$-module structure $\Afiltered\times \Mmit\rightarrow \Mmit$. 
Conversely, for any exhaustive left $\Afiltered$-module $\Mmit$, the action $\Afiltered\times \Mmit\rightarrow \Mmit$ restricts to an exhaustive left\/ $\Agraded$-module structure.
\end{proposition}
Thus the notion of exhaustive left $\Agraded$-modules and that of exhaustive left $\Afiltered$-modules have no essential differences.

\begin{note}
An exhaustive module is nothing else but a torsion module with respect to the left linear topology. (See \cite{Ga}.)
\end{note}
\subsection{Generalized eigenspaces of exhaustive modules}
\label{107}
Let $\Agraded$ be a compatible degreewise topological algebra.
For a left $\Agraded$-module $\Mmit$ and a nonnegative integer $n$, we set 
\begin{equation}
  \Krm_n(\Mmit)
  =\ensemble{v\in\Mmit}{$\Igraded_n(\Agraded)\cdot v=0$}.
\end{equation}
Note that $\Mmit$ is exhaustive if and only if $\Mmit=\bigcup_n\Krm_n(\Mmit)$ and if this is the case then 
$\Krm_n(\Mmit)
  =\ensemble{v\in\Mmit}{$\Fgraded_{-n-1}\Agraded\cdot v=0$}
$ by Lemma \ref{106}.

Consider the case when $(\Agraded,\ham)$ is a weakly quasi-finite algebra graded by Hamiltonian. 
Since the space $\Krm_n(\Mmit)$ has a structure of a left $\Qrm_n(0)$-module,  it decomposes into the sum of generalized eigenspaces: 
\begin{equation}
  \Krm_n(\Mmit)
  =\sum_{\lambda\in\Omega_n}\Krm_n(\Mmit)[\lambda]
.
\end{equation}
Hence if $\Mmit$ is exhaustive then the whole space $\Mmit$ also decomposes into the sum of generalized eigenspaces. 
Let us set
\begin{equation}
  \Erm_n(\Mmit)=\sum_{\lambda\in\Gamma_n}\Mmit[\lambda].
\end{equation}
The following is one of the key observations in the present paper.
\begin{proposition}
\label{050}
If $\Mmit$ is exhaustive then\/ $\Erm_n(\Mmit)\subset \Krm_n(\Mmit)\subset \Erm_{n+g}(\Mmit)$.
\begin{proof}
The containment $\Krm_n(\Mmit)\subset \Erm_{n+g}(\Mmit)$ follows from $\Omega_n\subset \Gamma_{n+g}$ because $\Krm_n(\Mmit)$ is a left $\Qrm_n(0)$-module.
The rest is similar to Lemma \ref{228}. 
Consider the space $\Mmit[\lambda]$ with $\lambda\in\Gamma_n$.
Then we have $\Agraded(d)\cdot \Mmit[\lambda]\subset \Mmit[\lambda+d]$.
Therefore, since $\lambda-n-1\notin \Gamma_\infty$, we have $\Fgraded_{-n-1}\Agraded\cdot \Mmit[\lambda]=0$.
This implies that $\Erm_n(\Mmit)\subset \Krm_n(\Mmit)$.
\end{proof}
\end{proposition}
\subsection{Equivalence of categories}
\label{121}
We will mean by the category of exhaustive left $\Agraded$-modules the full subcategory of the category of left $\Agraded$-modules consisting of exhaustive left $\Agraded$-modules.
Let $\Mmit$ be an exhaustive left $\Agraded$-module.
Then the space 
 $
  \Erm_n(\Mmit)
$ 
is a left $\Afinite_n$-module.
A homomorphism $\phi:\Mmit'\rightarrow \Mmit''$ of left $\Agraded$-modules induces a map
\begin{equation}
  \Erm_n(\phi):\Erm_n(\Mmit')\rightarrow \Erm_n(\Mmit'').
\end{equation}
Thus we have a functor $\Erm_n(-)$ from the category of exhaustive left $\Agraded$-module to the category of left $\Afinite_n$-modules.
\begin{lemma}
\label{097}
If $\Mmit$ is exhaustive then\/ $\Krm_0(\Mmit)=0$ implies $\Mmit=0$.
\begin{proof}
Suppose $\Krm_0(\Mmit)=0$ and $\Mmit\ne 0$. 
Let $v$ be a nonzero element of $\Mmit$.
Since $\Mmit$ is exhaustive, $\Igraded_{n}(\Agraded)\cdot v=0$ and hence $\Fgraded_{-n-1}\Agraded\cdot v=0$ for sufficiently large $n$. 
Take the maximal integer $n$ for which $\Fgraded_{-n-1}\Agraded\cdot v\ne 0$.
Then by the maximality $\Fgraded_{-1}\Agraded\cdot (\Fgraded_{-n-1}\Agraded\cdot v)=0$ and hence 
$\Igraded_0(\Agraded)\cdot (\Fgraded_{-n}\Agraded\cdot v)=0$ 
by Lemma \ref{106}.
Hence $\Fgraded_{-n}\Agraded\cdot v\subset \Krm_0(\Mmit)=0$ which is a contradiction.
\end{proof}
\end{lemma}
\begin{lemma}
\label{049}
A homomorphism\/ $\phi:\Mmit'\rightarrow \Mmit''$ of exhaustive left\/ $\Agraded$-modules is an isomorphism if and only if\/ $\Erm_\gap(\phi)$ is an isomorphism of vector spaces.
\begin{proof}
Assume that $\Erm_\gap(\phi)$ is an isomorphism of vector spaces.
Then, by Lemma \ref{097}, we see $\Ker \phi=0$ and $\Coker \phi=0$ since $\Krm_0(\Ker\phi)\subset \Erm_\gap(\Ker \phi)=\Ker(\Erm_\gap(\phi))=0$ and $\Krm_0(\Coker\phi)\subset \Erm_\gap(\Coker \phi)=\Coker(\Erm_\gap(\phi))=0$.
\end{proof}
\end{lemma}
Recall the $(\Agraded,\Afinite_n)$-bimodule $\Prm_n$ defined in Subsection \ref{100}, which is an exhaustive left $\Agraded$-module.
Therefore, for a left $\Afinite_n$-module ${X}$, the space 
 $
  \Prm_{n}\otimes_{\Afinite_n} {X}
$ 
is an exhaustive left $\Agraded$-module.
We note that $\Afinite_n=\Erm_n(\Prm_n)$. 
\begin{lemma}
\label{134}
For a left\/ $\Afinite_n$-module ${X}$, the map\/ 
 $
\Afinite_n\otimes_{\Afinite_n}{X}\rightarrow \Prm_{n}\otimes_{\Afinite_n}{X}
$ 
induced by the inclusion\/ $\Afinite_n\rightarrow \Prm_{n}$ is injective.
\begin{proof}
We set
 $
\Erm^\perp_n(\Prm_{n})=\sum_{\lambda\in\Gamma_\infty\setminus\Gamma_n}\Prm_{n}[\lambda]
$. 
Then the decomposition
 $
\Prm_{n}
=\Erm_n(\Prm_{n})\oplus \Erm^\perp_n(\Prm_{n})
=\Afinite_{n}\oplus \Erm^\perp_n(\Prm_{n})
$ 
is a direct sum decomposition of a right $\Afinite_n$-module. 
Hence the map
 $
\Erm_n(\Prm_{n})\otimes_{\Afinite_n}{X}\rightarrow \Prm_{n}\otimes_{\Afinite_n}{X}
$ 
is injective.
\end{proof}
\end{lemma}
Now we come to the second main result of Part I. 
Recall the number $\gap$ defined by (\ref{231}). 
\begin{theorem}
\label{146}
Let\/ $(\Agraded,\ham)$ be a weakly quasi-finite algebra graded by Hamiltonian and let $n$ be an integer such that $n\geq \gap$.
Then the functors $\Erm_n(-)$ and $\Prm_{n}\otimes_{\Afinite_n} -$ are mutually inverse equivalences of categories between the category of exhaustive left $\Agraded$-modules and the category of left $\Afinite_n$-modules.
\begin{proof}
Let ${X}$ be a left $\Afinite_n$-module. 
By Lemma \ref{134}, the map
 $
\Afinite_{n}\otimes_{\Afinite_n}{X}\rightarrow \Prm_{n}\otimes_{\Afinite_n}{X}
$ 
is injective.
Therefore
\begin{equation}
\label{234}
  {X}
  \simeq\Afinite_n\otimes_{\Afinite_n}{X}
  \simeq \Erm_n(\Prm_{n}\otimes_{\Afinite_n}{X})
.
\end{equation}
Let $\Mmit$ be an exhaustive left $\Agraded$-module. Then by letting ${X}=\Erm_n(\Mmit)$ in (\ref{234}), we have
\begin{equation}
  \Erm_n(\Mmit)
  \simeq \Erm_n(\Prm_{n}\otimes_{\Afinite_n}\Erm_n(\Mmit))
.
\end{equation}
Since $n\geq \gap$, we have
 $
  \Prm_{n}\otimes_{\Afinite_n}\Erm_n(\Mmit)
  \simeq\Mmit
$ by Lemma \ref{049}.
\end{proof}
\end{theorem}
\begin{corollary}
\label{191}
Let\/ $(\Agraded,\ham)$ be a weakly quasi-finite algebra graded by Hamiltonian and let $n$ be an integer such that $n\geq \gap$.
Then the module\/ $\Prm_n$ is a progenerator of the category of exhaustive left\/ $\Agraded$-modules.
\end{corollary}
Let us consider the case when $\Agraded$ is quasi-finite. 
Then the algebra $\Agraded_n$ is finite-dimensional for all $n$ by Proposition  \ref{232}. 
Thus we have the following corollary.
\begin{corollary}
\label{223}
If\/ $(\Agraded,\ham)$ is a quasi-finite algebra graded by Hamiltonian then the category of exhaustive left\/ $\Agraded$-modules is equivalent to the category of left modules over a finite-dimensional algebra. 
\end{corollary}
\vskip2ex
We may likewise define the notion of exhaustive right $\Agraded$-modules and exhaustive right $\Afilteredv$-modules by using the spaces $\Igraded^\vee_n(\Agraded)$ and $\Ifilteredv_n(\Afilteredv)$, respectively.
We have analogous results for the right modules as well.
\begin{note}
Theorem \ref{146} may be understood to be a particular case of a topological variant of Morita equivalences. 
See \cite{Ga} for general theory of equivalences between abelian categories and module categories and \cite{Gr} and \cite{Mez} for results closely related to Theorem \ref{146}. 
\end{note}
%
%
\section{Coexhaustive modules}
\label{054}
We will now consider the notion of coexhaustive modules, which is {\it dual\/} to that of exhaustive modules. 
In this section, we will give the definition of coexhaustive modules and describe their topologies by means of the generalized eigenspaces. 
The precise duality will be considered in the next section under necessary finiteness assumptions. 
\subsection{Coexhaustive modules}
\label{055}
Let $(\Agraded,\ham)$ be a weakly quasi-finite algebra graded by Hamiltonian and let $\Afilteredv$ be the right filterwise completion. 
Let $\Mcal$ be a left $\Afilteredv$-module endowed with a linear topology such that the action $\Afilteredv\times\Mcal\rightarrow \Mcal$ is continuous. 
Let us denote by $\Ifilteredv_n(\Mcal)$ the closure of the space $ \Ffilteredv_{-n-1}\Afilteredv\cdot\Mcal$.
We assume that the sequence $\{\Ifilteredv_n(\Mcal)\}$ forms a fundamental system of open neighborhoods of zero. 
\begin{definition}
A {\it  compatible topological left $\Afilteredv$-module\/} is a left $\Afilteredv$-module endowed with a linear topology which satisfies the conditions mentioned above.
\end{definition}
Let $\Mcal$ be a compatible topological left $\Afilteredv$-module.
Then we have
\begin{equation}
  \Ffilteredv_p\Afilteredv\cdot \Ifilteredv_n(\Mcal)
    \subset \Ifilteredv_{n-p}(\Mcal)
,\quad
  \Ifilteredv_n(\Afilteredv)\cdot \Mcal
    \subset \Ifilteredv_n(\Mcal)
\end{equation}
because the action $\Afilteredv\times\Mcal\rightarrow \Mcal$ is continuous. 
\begin{definition}
A compatible topological left $\Afilteredv$-module $\Mcal$ is {\it  coexhaustive\/} if $\Mcal$ is complete as a topological vector space. 
\end{definition}

Obviously, the completion of a compatible topological left $\Afilteredv$-module has a canonical structure of a coexhaustive left $\Afilteredv$-module. 
\subsection{Generalized eigenspaces of coexhaustive modules}
\label{018}
Let $\Mcal$ be a  coexhaustive left $\Afilteredv$-module.
Consider the space 
\begin{equation}
  \Qrm^\vee_n(\Mcal)=\Mcal/\Ifilteredv_n(\Mcal)
.
\end{equation}
Then this is a left $\Qrm^\vee_n(0)$-module and hence has a generalized eigenspace decomposition of the form
\begin{equation}
  \Qrm^\vee_n(\Mcal)
    =\sum_{\mu\in\Omega_n}\Qrm^\vee_n(\Mcal)[\mu].
\end{equation}
\begin{lemma}
\label{080}
If\/ $\lambda\in \Gamma_m$ and $n\geq m$ then the restriction\/ 
 $ 
  \Qrm^\vee_n[\lambda]
    \rightarrow \Qrm^\vee_m[\lambda]
 $ 
of the canonical surjection is an isomorphism.
\begin{proof}
See the proof of Lemma \ref{074}. 
\end{proof}
\end{lemma}
This lemma implies that the subspace $\sum_{\lambda\in\Gamma_\infty}\Mcal[\lambda]$ is dense in $\Mcal$.
It also implies that the space 
\begin{equation}
  \Erm_n(\Mcal)=\sum_{\lambda\in\Gamma_n}\Mcal[\lambda]
\end{equation}
is a discrete subspace for each $n$.
Consider the space $\Erm^\perp_n(\Mcal)=\sum_{\lambda\in\Gamma_\infty\setminus\Gamma_n}\Mcal[\lambda]$ and let $\widehat{\Erm^\perp_n(\Mcal)}$ be its closure in $\Mcal$. 
Then the quotient 
\begin{equation}
  \Prm_n(\Mcal)=\Mcal/\widehat{\Erm^\perp_n(\Mcal)}
\end{equation}
is canonically isomorphic to $\Erm_n(\Mcal)$.
\begin{lemma}
\label{079}
If\/ $\Mcal$ is a  coexhaustive left\/ $\Afilteredv$-module then the canonical map\/ 
$\Erm_m(\Mcal)\rightarrow \Qrm^\vee_m(\Mcal)$ is injective and the canonical map\/ $\Qrm^\vee_m(\Mcal)\leftarrow \Erm_{m+g}(\Mcal)$ 
is surjective.
\begin{proof}
By Lemma \ref{080}, we know that the map $\Mcal[\lambda]\rightarrow \Qrm^\vee_n(\Mcal)[\lambda]$ is injective if $\lambda\in\Gamma_n$.
Hence the map $\Erm_n(\Mcal)\rightarrow \Qrm^\vee_n(\Mcal)$ is injective.
Now recall that $\Qrm^\vee_n(\Mcal)=\sum_{\lambda\in \Gamma_{n+\gap}}\Qrm^\vee_n(\Mcal)[\lambda]$.  
Hence the map $\Qrm^\vee_m(\Mcal)\leftarrow \Erm_{m+g}(\Mcal)$ is surjective.
\end{proof}
\end{lemma}
\subsection{Exhaustion and coexhaustion}
\label{081}
We will mean by the category of coexhaustive left $\Afilteredv$-modules the category for which the objects are the coexhaustive left $\Afilteredv$-modules and the morphisms are the continuous homomorphisms of left $\Afilteredv$-modules. 

\vskip2ex

Let $\Mmit$ be an exhaustive left $\Afiltered$-module and regard it as a left $\Agraded$-module by identifying $\Agraded$ with a subalgebra of $\Afiltered$ via the canonical map $\Agraded\rightarrow \Afiltered$.
We give $\Mmit$ the linear topology defined by  $\Igraded^\vee_n(\Mmit)=\Fgraded^\vee_{-n-1}\Agraded\cdot \Mmit$. 
Then $\Mmit$ becomes a left $\Agraded$-module such that the action $\Agraded\times \Mmit\rightarrow \Mmit$ is continuous. 
We let $\Qscr^\vee_\infty(\Mmit)$ be the completion: 
\begin{equation}
  \Qscr^\vee_\infty(\Mmit)=\varprojlim_n \Qrm^\vee_n(\Mmit)
,\quad
  \Qrm^\vee_n(\Mmit)=\Mmit/\Igraded^\vee_n(\Mmit)
.
\end{equation}
Then $\Qscr^\vee_\infty(\Mmit)$ is a coexhaustive left $\Afilteredv$-module. 

Let $\Mmit'$ and $\Mmit''$ be exhaustive left $\Afiltered$-modules and let $\phi:\Mmit'\rightarrow \Mmit''$ be a homomorphism of left $\Afiltered$-modules.
Then $\phi$ gives rise to a continuous homomorphism of left $\Agraded$-modules. 
Hence it induces a continuous homomorphism $\Qscr^\vee_\infty(\phi):\Qscr^\vee_\infty(\Mmit')\rightarrow \Qscr^\vee_\infty(\Mmit'')$ of coexhaustive left $\Afilteredv$-modules. 
We call the functor $\Qscr^\vee_\infty(-)$ the {\it coexhaustion functor}.
\vskip2ex
Conversely, let $\Mcal$ be a coexhaustive left $\Afilteredv$-module and regard it as a left $\Agraded$-module via the canonical map $\Agraded\rightarrow \Afilteredv$.
Consider the space
\begin{equation}
  \Krm_\infty(\Mcal)=\bigcup_n\Krm_n(\Mcal)
,\quad
  \Krm_n(\Mcal)
  =\ensemble{v\in\Mcal}{$\Igraded_n(\Agraded)\cdot v=0$}
.
\end{equation}
Then this is an exhaustive left $\Agraded$-module and hence an exhaustive left $\Afiltered$-module. 
Let $\varphi:\Mcal'\rightarrow \Mcal''$ be a continuous homomorphism of coexhaustive left $\Afilteredv$-modules. 
We regard $\varphi$ as a homomorphism of left $\Agraded$-modules via the canonical map $\Agraded\rightarrow \Afilteredv$. 
If $\Igraded_n(\Agraded)\cdot v=0$ then $\Igraded_n(\Agraded)\cdot \varphi(v)=\varphi(\Igraded_n(\Agraded)\cdot v)=\varphi(0)=0$. 
Therefore, $\varphi$ restricts to a homomorphism $\Krm_\infty(\varphi):\Krm_\infty(\Mcal')\rightarrow \Krm_\infty(\Mcal'')$ of left $\Agraded$-modules and hence a homomorphism of left $\Afiltered$-modules.
We call the functor $\Krm_\infty(-)$ the {\it exhaustion functor}.

\vskip2ex

Let us set
\begin{equation}
  \Pscr_\infty(\Mmit)=\varprojlim_n \Prm_n(\Mmit)
,\quad
  \Erm_\infty(\Mcal)=\bigcup_n \Erm_n(\Mcal)
.
\end{equation}
\begin{lemma}
\label{128}
For an exhaustive left\/ $\Afiltered$-module\/ $\Mmit$ its coexhaustion\/ $\Qscr^\vee_\infty(\Mmit)$ is canonically isomorphic to\/ $\Pscr_\infty(\Mmit)$  as topological vector spaces. 
For a coexhaustive left $\Afilteredv$-module $\Mcal$ its exhaustion $\Krm_\infty(\Mcal)$ is canonically isomorphic to $\Erm_\infty(\Mcal)$ as vector spaces. 
\begin{proof}
Conclusions are clear by the arguments in Subsections \ref{107} and \ref{018}.
\end{proof}
\end{lemma}
Now the following theorem is an immediate consequence of this lemma.
\begin{theorem}
\label{084}
Let\/ $(\Agraded,\ham)$ be a weakly quasi-finite algebra graded by Hamiltonian 
 and let $\Afiltered$ and $\Afilteredv$ be the left and the right filterwise completions, respectively. 
Then the functors $\Qscr^\vee_\infty(-)$ and\/ $\Krm_\infty(-)$ are mutually inverse equivalences of categories between the category of exhaustive left $\Afiltered$-modules and the category of  coexhaustive left $\Afilteredv$-modules.
\end{theorem}
In particular, the category of coexhaustive left $\Afilteredv$-modules is an abelian category. 
\vskip2ex
We may likewise define the notion of  coexhaustive right $\Agraded$-modules and  coexhaustive right $\Afiltered$-modules by using the spaces $\Igraded_n(\Agraded)$ and $\Ifiltered_n(\Afiltered)$, respectively.
We have analogous results for the right modules as well.
%
%
\section{Duality for quasi-finite modules}
\label{115}
We now consider finiteness conditions for exhaustive modules and for  coexhaustive modules and discuss the duality between the categories of such modules.
\subsection{Quasi-finiteness of exhaustive modules}
\label{131}
Let $(\Agraded,\ham)$ be a weakly quasi-finite algebra graded by Hamiltonian and let $\Afiltered$ and $\Afilteredv$ be the left and the right filterwise completions, respectively. 
\begin{definition}
A {\it quasi-finite left $\Afiltered$-module\/} is an exhaustive left $\Afiltered$-module $\Mmit$ such that the spaces $\Krm_n(\Mmit)$ are finite-dimensional for all $n$. 
\end{definition}
For an exhaustive right $\Afilteredv$-module $\Nmit$, we define the space $\Krm^\vee_n(\Nmit)$ in the same way as $\Krm_n(\Mmit)$ for an exhaustive left $\Afiltered$-module $\Mmit$:
\begin{equation}
  \Krm^\vee_n(\Nmit)
    =\ensemble{v\in\Nmit}{$v\cdot \Ifilteredv_n(\Afilteredv)=0$}
.
\end{equation}
\begin{definition}
A {\it quasi-finite right $\Afilteredv$-module\/} is an exhaustive right $\Afilteredv$-module $\Nmit$ such that the spaces $\Krm^\vee_n(\Nmit)$ are finite-dimensional for all $n$. 
\end{definition}

The following proposition characterizes the quasi-finiteness of an algebra by means of the quasi-finiteness of the canonical quotient modules.  

\begin{proposition}
\label{152}
Let\/ $(\Agraded,\ham)$ be a weakly quasi-finite algebra with Hamiltonian. 
Then the following conditions are equivalent:
\begin{enumerate}[\rm (a)]
\item
$\Agraded$ is a quasi-finite algebra.
\item
The left canonical quotient modules\/ $\Qrm_n$ are quasi-finite for all $n$. 
\item
The right canonical quotient modules\/ $\Qrm^\vee_n$ are quasi-finite for all $n$. 
\end{enumerate}
\begin{proof}
We will show the equivalence of (a) and (b). 
Since $\Krm_m(\Qrm_n)$ is an $(\Hrm_m,\Hrm_n)$-bimodule, the pair of the left and the right eigenvalues of $\ham$ on $\Krm_m(\Qrm_n)$ are in the finite set $\Omega_m\times \Omega_n$.
Then if $\Krm_m(\Qrm_n)(d)\ne 0$ then $d=\lambda-\mu$ for some $\lambda\in \Omega_m$ and $\mu\in\Omega_n$.
Therefore,  we have
\begin{equation}
  \Krm_m(\Qrm_n)\subset\bigoplus_{-n\leq d\leq m+\gap}\Qrm_n(d).
\end{equation}
On the other hand, it is easy to see that 
 $
  \Qrm_n(d)\subset \Krm_{n+d+1}(\Qrm_n).
$ 
Therefore we immediately see that  $\Krm_m(\Qrm_n)$ are finite-dimensional for all $m$ and $n$ if and only if $\Qrm_n(d)$ are  finite-dimensional for all $d$ and $n$. The proof for the equivalence of (a) and (c) is similar.

\end{proof}
\end{proposition}
Now consider the case when $(\Agraded,\ham)$ is a quasi-finite algebra graded by Hamiltonian. 
Since $\Erm_n(\Mmit)$ is a left $\Afinite_n$-module and $\Afinite_n$ is 
finite-dimensional, $\Erm_n(\Mmit)$ is finite-dimensional if and only if it 
is finitely generated as a left $\Afinite_n$-module. 
\begin{proposition}
\label{150}
Let\/ $(\Agraded,\ham)$ be a quasi-finite algebra graded by Hamiltonian.
Then the following conditions for an exhaustive left $\Afiltered$-module $\Mmit$ are equivalent: 
\begin{enumerate}[\rm(a)]
\topsep=5ex
\item
$\Mmit$ is a quasi-finite left $\Afiltered$-module.
\item
$\Mmit$ is finitely generated as a left $\Afiltered$-module.
\item
$\Erm_n(\Mmit)$ are finite-dimensional for all $n$. 
\item
$\Erm_\gap(\Mmit)$ is finite-dimensional.
\end{enumerate}
\begin{proof}
We will show (a)$\Rightarrow$(c)$\Rightarrow$(d)$\Rightarrow$(b)$\Rightarrow$(a).
Assume that $\Mmit$ is quasi-finite.
Then $\Erm_n(\Mmit)$ are finite-dimensional since $\Erm_n(\Mmit)\subset \Krm_n(\Mmit)$ by Proposition \ref{050}. 
In particular, $\Erm_\gap(\Mmit)$ is finite-dimensional. 
Next assume that $\Erm_\gap(\Mmit)$ is finite-dimensional and let $\Mmit'$ be the left $\Afiltered$-submodule of $\Mmit$ generated by $\Erm_\gap(\Mmit)$.
Then $\Krm_0(\Mmit/\Mmit')\subset \Erm_\gap(\Mmit/\Mmit')=0$. 
By Lemma \ref{097}, we have $\Mmit=\Mmit'$ and hence $\Mmit$ is finitely generated.
Now assume that $\Mmit$ is finitely generated and let $v_1,\ldots,v_k$ be a set of generators. 
Since $\Mmit$ is exhaustive, there exists $n_1,\ldots,n_k$ such that $\Igraded_{n_i}(\Afiltered)\cdot v_i=0$ for $i=1,\ldots,k$.
This implies the existence of a surjective homomorphism 
 $ 
  \Qrm_{n_1}\times\cdots\times \Qrm_{n_k}\rightarrow \Mmit
$  
of left $\Afiltered$-modules. 
Since $\Agraded$ is quasi-finite and hence the modules $\Qrm_{n_i}$ are quasi-finite, so is the image $\Mmit$.
\end{proof}
\end{proposition}
\subsection{Quasi-finiteness of coexhaustive modules}
Let us now turn to the quasi-finiteness of coexhaustive modules. 
\begin{definition}
A {\it quasi-finite left $\Afilteredv$-module\/} is a coexhaustive left $\Afilteredv$-module $\Mcal$ such that the spaces $\Qrm^\vee_n(\Mcal)$ are finite-dimensional for all $n$. 
\end{definition}

For a coexhaustive right $\Afiltered$-module $\Ncal$, we set
\begin{equation}
  \Qrm_n(\Ncal)=\Ncal/\Ifiltered_n(\Ncal)
.
\end{equation}
\begin{definition}
A {\it quasi-finite right $\Afiltered$-module\/} is a coexhaustive right $\Afiltered$-\break module $\Ncal$ such that the spaces $\Qrm_n(\Ncal)$ are finite-dimensional for all $n$. 
\end{definition}
For an exhaustive left $\Afiltered$-module $\Mmit$, the canonical map $\Erm_n(\Mmit)\rightarrow \Prm^\vee_n(\Mmit)$ is an isomorphism. 
Therefore, Lemma \ref{079} implies that $\Erm_n(\Mmit)$ are finite-dimensional for all $n$ if and only if $\Prm^\vee_n(\Mmit)$ are finite-dimensional for all $n$.
Therefore, the consideration in Subsection \ref{081} implies the following result.
\begin{proposition}
\label{151}
Let\/ $(\Agraded,\ham)$ be a quasi-finite algebra graded by Hamiltonian.
Let $\Mmit$ be an exhaustive left $\Afiltered$-module and let $\Mcal$ be its coexhaustion.
Then the exhaustive module $\Mmit$ is  quasi-finite if and only if the  coexhaustive module $\Mcal$ is quasi-finite.
\end{proposition}
We have analogous results for quasi-finite right $\Afilteredv$-modules and for quasi-finite right $\Afiltered$-modules.
\subsection{Duality}
Let $(\Agraded,\ham)$ be a weakly quasi-finite algebra graded by Hamiltonian. 
Let $\Mmit$ be an exhaustive left $\Afiltered$-module and consider its full dual space $\Ncal$: 
\begin{equation}
  \Ncal=\Mmit^*=\Hom_{\KK}(\Mmit,\KK).
\end{equation}
Then this space becomes a right $\Afiltered$-module. 

Since $\Mmit$ is exhaustive, it has the generalized eigenspace decomposition 
 $
\Mmit=\bigoplus_{\lambda\in\Gamma_\infty}\Mmit[\lambda]
$, 
which gives rise to 
\begin{equation}
  \Ncal=\prod_{\lambda\in\Gamma_\infty}\Mmit[\lambda]^*
\end{equation}
where
 $\Mmit[\lambda]^*$ is the set of linear functions $f:\Mmit\rightarrow \KK$ such that $f(\Mmit[\mu])=0$ holds for any $\mu\ne\lambda$.
Therefore, by setting 
\begin{equation}
  \Jrm_n(\Ncal)
    =\ensemble{f\in\Ncal}{$f(\Erm_n(\Mmit))=0$}
,
\end{equation}
we have 
\begin{equation}
  \Ncal/\Jrm_n(\Ncal)
    =\prod_{\lambda\in\Gamma_n}\Mmit[\lambda]^*
.
\end{equation}
We give $\Ncal$ the linear topology defined by $\Jrm_n(\Ncal)$. 

\begin{lemma}
\label{155}
The topological right\/ $\Afiltered$-module\/ $\Ncal$ is a coexhaustive right\/ $\Afiltered$-module.
\begin{proof}
Completeness of $\Ncal$ as a topological vector space is clear from the definition. 
Since 
 $
  (f\cdot a)(\Erm_n(\Mmit))=f(a\cdot \Erm_n(\Mmit))=0
$ 
holds for any $f\in\Ncal$ and any $a\in\Ifiltered_n(\Afiltered)$, we have
 $
  \Ncal\cdot\Ifiltered_n(\Afiltered)
    \subset \Jrm_n(\Ncal)
$.
Hence the action $\Ncal\times \Afiltered\rightarrow \Ncal$ is continuous.
Let $\Ifiltered_n(\Ncal)$ be the closure of $\Ncal\cdot\Ifiltered_n(\Afiltered)$.
It remains to show that $\{\Ifiltered_n(\Ncal)\}$ forms a fundamental system of open neighborhoods of zero. 
Since $\Jrm_n(\Ncal)$ is closed, we have $\Ifiltered_n(\Ncal)\subset \Jrm_n(\Ncal)$. 
Consider the quotient $\Ncal/\Ifiltered_n(\Ncal)$.
Since we have 
 $
  \Ifiltered_n(\Ncal)\cdot \Ifiltered_n(\Afiltered)\subset \Ifiltered_n(\Ncal)
$, 
the quotient space $\Ncal/\Ifiltered_n(\Ncal)$ is a right $\Qrm_n(0)$-module.
In particular, it is a right $\Hrm_n$-module.
Hence it has a generalized eigenspace decomposition with respect to the right action of $\ham$ with the eigenvalues 
belonging to $\Omega_n\subset \Gamma_{n+\gap}$. 
Since $\Jrm_{n+\gap}(\Ncal)$ is the closure of $\Erm^\perp_{n+\gap}(\Ncal)$, 
we have $\Jrm_{n+\gap}(\Ncal)\subset \Ifiltered_n(\Ncal)$. 
Hence $\Ifiltered_n(\Ncal)$ with $n\in\NN$ form a basis of open neighborhoods of zero.
\end{proof}
\end{lemma}
Thus we have shown that for any exhaustive left $\Afiltered$-module, its full dual naturally becomes a  coexhaustive right $\Afiltered$-module.
\vskip2ex
Conversely, let $\Ncal$ be a  coexhaustive right $\Afiltered$-module and consider the continuous dual space $\Mmit$:
\begin{equation}
  \Mmit=\Hom^{\rm cont}_{\KK}(\Ncal,\KK).
\end{equation}
Recall that the base field $\KK$ is given the discrete topology.

Take any element of $\Mmit$, which is given by a continuous map
 $
  f:\Ncal\rightarrow \KK.
$ 
Then, for any $a\in\Afiltered$, the map
 $
  a\cdot f:\Ncal\rightarrow \KK
$ 
is also continuous since, for $a\in\Fgraded_p\Afiltered$, we have
\begin{equation}
  (a\cdot f)(v+\Irm_n(\Ncal))
    \subset f(v\cdot a)+f(\Irm_{n-p}(\Ncal))
    =f(v\cdot a)
\end{equation}
for sufficiently large $n$.
We then have 
\begin{equation}
  (\Igraded_{n}(\Afiltered)\cdot f)(v)
  =f(v\cdot \Igraded_{n}(\Afiltered))
  \subset f(\Ifiltered_{n}(\Ncal))
  =0
\end{equation}
for all $v$. 
Therefore, the module $\Mmit$ is exhaustive. 

\vskip2ex

Thus we have defined contravariant functors $\Hom_\KK(-,\KK)$ and $\Hom^{\rm cont}_{\KK}(-,\KK)$ between the category of exhaustive left $\Afiltered$-modules and the category of coexhaustive right $\Afiltered$-modules. 
\begin{proposition}
Let\/ $(\Agraded,\ham)$ be a weakly quasi-finite algebra graded by Hamiltonian. 
Then the functors\/ $\Hom_\KK(-,\KK)$ and\/ $\Hom^{\rm cont}_{\KK}(-,\KK)$ give rise to mutually inverse duality of categories between the category of quasi-finite left $\Afiltered$-modules and the category of quasi-finite right $\Afiltered$-modules.
\end{proposition}
\begin{note}
\label{164}
A quasi-finite right $\Afiltered$-module is nothing else but a linearly compact right $\Afiltered$-module. 
We may understand the above-mentioned duality as a version of the Lefschetz duality \cite{Le} between discrete modules and linearly compact modules. 
See also \cite{Mac}. 
\end{note}
\subsection{Involution}
\label{031}
Let $(\Agraded,\ham)$ be a quasi-finite algebra graded by Hamiltonian.
Suppose given a linear involution
 $
  \theta:\Agraded\rightarrow\Agraded
$ 
satisfying the following conditions:
\begin{enumerate}[\rm(i)]
\item
$\theta(a\cdot b)=\theta(b)\cdot \theta(a)$ for any $a,b\in \Agraded$ and $\theta(\ham)=\ham$.
\item
$\theta:\Agraded(d)\rightarrow\Agraded(-d)$ is continuous.
\end{enumerate}
Note that by (i) we have $[\ham,\theta(a)]=[\theta(h),\theta(a)]=-\theta([\ham,a])$.
Hence $\theta(\Agraded(d))=\Agraded(-d)$ and the condition (ii) makes sense. 
We then have that $\theta(\Igraded_n(\Agraded))=\Igraded^\vee_n(\Agraded)$.
Hence $\theta$ extends to anti-isomorphisms
\begin{equation}
  \theta_\infty:\Afiltered\rightarrow \Afilteredv
,\quad
  \theta^\vee_\infty:\Afilteredv\rightarrow \Afiltered
\end{equation}
of filtered topological algebras such that $\theta_\infty\compo\theta^\vee_\infty=1$ and $\theta^\vee_\infty\compo\theta_\infty=1$.
We will denote the maps $\theta_\infty$ and $\theta^\vee_\infty$ by the same symbol $\theta$ by abuse of notation.
Let $\Mmit$ be a left $\Afiltered$-module and let $\vartheta(\Mmit)$  be the same space $\Mmit$ as a vector space. 
We give $\vartheta(\Mmit)$ a structure of a right $\Afilteredv$-module by letting $v\cdot a=\theta(a)\cdot v$ for $a\in\Afilteredv$ and $v\in\Mmit$. 
Similarly, for a right $\Afilteredv$-module $\Nmit$, we define a left $\Agraded$-module $\vartheta^\vee(\Nmit)$ in a similar way.
\begin{proposition}
The functors $\vartheta$ and $\vartheta^\vee$ are mutually inverse equivalences of categories between the category of left $\Afiltered$-modules and the category of right $\Afilteredv$-modules.
\end{proposition}
Let us compose the involution, the duality and the exhaustion. 
Then we get an auto-duality 
\begin{equation}
  \Drm:\Mmit\mapsto \vartheta(\Krm_\infty(\Hom_\KK(\Mmit,\KK)))
\end{equation}
of the category of quasi-finite left $\Afiltered$-modules. 
Note that $\Krm_\infty(\Hom_\KK(\Mmit,\KK))$ is the restricted dual space: 
\begin{equation}
  \Krm_\infty(\Hom_\KK(\Mmit,\KK))
    =\ensemble{f}{$f(\Erm^\perp_n(\Mmit))=0$ for some $n$}.
\end{equation}
\subsection{Finiteness theorems}
\label{224}
Let $(\Agraded,\ham)$ be a quasi-finite algebra and consider the following categories. 
\begin{enumerate}[\rm L1.]
\item
The category of quasi-finite left $\Afiltered$-modules.
\item
The category of quasi-finite left $\Afilteredv$-modules.
\item
The category of finitely generated left $\Afinite_n$-modules.
\end{enumerate}
\begin{enumerate}[\ \rm R1.]
\item
The category of quasi-finite right $\Afilteredv$-modules.
\item
The category of quasi-finite right $\Afiltered$-modules.
\item
The category of finitely generated right $\Afinite_n$-modules.
\end{enumerate}
Recall that quasi-finite left $\Afiltered$-modules and 
 right $\Afilteredv$-modules are exhaustive whereas  quasi-finite left
 $\Afilteredv$-modules and right $\Afiltered$-modules are coexhaustive. 
\vskip2ex
The following theorem summarizes the results obtained so far regarding the quasi-finite modules over quasi-finite algebras.
\begin{theorem}
\label{233}
Let\/ $(\Agraded,\ham)$ be a quasi-finite algebra with Hamiltonian and let $n$ be an integer such that $n\geq \gap$.
\begin{enumerate}[\rm(1)]
\item
The categories {\rm L1}, {\rm L2} and {\/\rm L3} are equivalent to each other. 
\item
The categories {\rm R1}, {\rm R2} and {\/\rm R3} are equivalent to each other. 
\item
The categories {\rm L1}, {\rm L2} and {\/\rm L3} and the categories {\rm R1}, {\rm R2} and {\/\rm R3} are dual to each other. 
\item
If $\Agraded$ has an involution $\theta$ then the categories {\rm L1}, {\rm L2} and {\/\rm L3} and the categories {\rm R1}, {\rm R2} and {\/\rm R3} are equivalent to each other. 
\end{enumerate}
\end{theorem}
%
\part{Quasi-finiteness and Zhu's finiteness condition}
%
%
%
\section{Vertex operator algebras and current algebras}
\label{028}
We now turn our attention to vertex operator algebras. 
In this section, we will describe in detail the construction and properties of the universal enveloping algebra associated with a vertex operator algebra, which we will simply call the {\it current algebra}, in order to give precise statements which seem to have been overlooked in the literatures. 
\subsection{Vertex operator algebra}
\label{168}
Recall that a vertex operator algebra is a graded vector space $\VOA$ with the 
grading being indexed by integers equipped with countably many bilinear maps 
indexed by integers and two distinguished elements, called the vacuum vector 
and the conformal vector (or the Virasoro element), satisfying a number of 
axioms (\cite{Bor}, \cite{FLM}, \cite{FHL}), which we will describe briefly 
below. 
See \cite{MN} for an account. 

Let us denote the homogeneous subspaces of the grading of $\VOA$ by $\VOA[r]$ with $r\in\ZZ$. 
It is assumed that there exists an $m\in\NN$ such that $\VOA[r]=0$ for $r<-m$. 
Therefore, the grading of $\VOA$ is of the following form
\begin{equation}
  \VOA=\bigoplus_{r=-m}^\infty\VOA[r]
.
\end{equation}
We will write $\Delta({u})=r$ when ${u}$ belongs to the subspace $\VOA[r]$ and call $\Delta({u})$ the {\it weight\/} of the element ${u}$. 

Let us denote the countably many bilinear maps by 
\begin{equation}
\label{170}
  \mu_n:\VOA\times \VOA\rightarrow \VOA
,\quad 
  ({u},{v})\mapsto {u}_{(n)}{v}
.
\end{equation}
It is assumed that they satisfy $\VOA[r]_{(n)}\VOA[s]\subset \VOA[r+s-n-1]$.
In other words, for homogeneous $u$ and $v$, we have
\begin{equation}
\Delta(u_{(n)}v)=\Delta(u)+\Delta(v)-n-1
.
\end{equation}
Then the sums in the following expression are finite:
\begin{equation}
\label{169}
\eqalign{
%
%
   \sum_{i=0}^\infty \binom{p}{i}({u}_{(r+i)}{v})_{(p+q-i)}{w}
%
%
&
\cr&\hskip-6em
  =\sum_{i=0}^\infty (-1)^i\binom{r}{i}\bigl({u}_{(p+r-i)}({v}_{(q+i)}{w})
%
%
  -(-1)^r{v}_{(q+r-i)}({u}_{(p+i)}{w})\bigr).
\cr}
\end{equation}
The bilinear maps (\ref{170}) are assumed to satisfy this identity for any ${u},{v},{w}\in \VOA$ and any $p,q,r\in\ZZ$. 
This set of identities is the main identity of vertex operator algebras called the {\it Borcherds identity\/} or the {\it Cauchy-Jacobi\/} identity.
The vacuum vector $\vac$ is an element of weight $0$. 
It enjoys, for any ${u}\in\VOA$, the relation ${u}_{(-1)}\vac={u}$ and ${u}_{(n)}\vac=0$ for $n\geq 0$.
We set
\begin{equation}
  \Der {u}={u}_{(-2)}\vac.
\end{equation}
Then it follows that the operator $\Der:\VOA\rightarrow \VOA$ is a derivation with respect to the operations $_{(n)}$ for every $n$.
The conformal vector $\omega$ is an element of weight $2$. 
It satisfies
\begin{equation}
  \omega_{(n)}\omega=0,\ (n\geq 4 \hbox{ or } n=2)
,\ 
  \omega_{(1)}\omega=2\omega
,\ 
  \omega_{(3)}\omega\in\KK \vac
.
\end{equation}
Then it follows that the operators $\Lvir_n:\VOA\rightarrow \VOA$ defined by $\Lvir_n{u}=\omega_{(n+1)}{u}$ satisfy the Virasoro commutation relation with  the central charge $\crm$ given by $2\omega_{(3)}\omega=\crm\,\vac$. 
Among the Virasoro operators $\Lvir_n$, the $\Lvir_0$ and $\Lvir_{-1}$ have special roles: 
the weights of $\VOA$ are given by the eigenvalues of $\Lvir_0$ and the derivation $\Der$ agrees with $\Lvir_{-1}$. 
In other words, 
\begin{equation}
  \VOA[r]=\ensemble{v\in\VOA}{$\Lvir_0v=rv$}
,\quad
  \Der{u}=\Lvir_{-1}{u}.
\end{equation}
The weight subspaces are usually assumed to be finite-dimensional. 
Among the conditions regarding the Virasoro operators, only the one concerning  the eigenvalues of $\Lvir_0$ will be used in the rest of the paper. 
\subsection{The current Lie algebras}
\label{029}
Let us consider the space
\begin{equation}
  \VOA[t,t^{-1}]=\VOA\otimes_{\KK}\KK[t,t^{-1}].
\end{equation}
We define a bilinear map $\VOA[t,t^{-1}]\times\VOA[t,t^{-1}]\rightarrow \VOA[t,t^{-1}]$ by setting
\begin{equation}
\label{194}
  [{u}\otimes t^m,{v}\otimes t^n]=\sum_{i=0}^\infty\binom{m}{n}({u}_{(i)}{v})\otimes t^{m+n-i-1}.
\end{equation}
%
%
%
%
%
Consider the quotient space
\begin{equation}
  \ggm=\VOA[t,t^{-1}]\,/\,\partial\VOA[t,t^{-1}]
\end{equation}
where $\partial:\VOA[t,t^{-1}]\rightarrow \VOA[t,t^{-1}]$ is defined by
\begin{equation}
  \partial({u}\otimes t^n)=\Der {u}\otimes t^n+n{u}\otimes t^{n-1}.
\end{equation}
Then the bracket operation on $\VOA[t,t^{-1}]$ induces a bilinear operation on $\ggm$ denoted by the same symbol which gives a structure of a Lie algebra on $\ggm$ (\cite{Bor}).
We will call the Lie algebra $\ggm$ the associated {\it current Lie algebra}.
We will denote the image of an element of $\VOA[t,t^{-1}]$ in $\ggm$ by the same symbol. 
Since $n\,\vac\otimes t^{n+1}=\partial(\vac\otimes t^n)$, we know that $\vac\otimes t^n=0$ in $\ggm$ is zero unless $n=-1$, when $\vac\otimes t^{-1}$ is central.
It will be useful to introduce the following notation:
\begin{equation}
  J_n({u})={u}\otimes t^{n+\Delta({u})-1}
\end{equation}
for a homogeneous ${u}$ and extend it linearly. 
We denote its image in $\ggm$ by the same symbol and assign the degree $-n$ to $J_n({u})$. 
Then the associated current Lie algebra is graded by the degree:
\begin{equation}
  \ggm=\bigoplus_{d\in\ZZ}\ggm(d).
\end{equation}
Note the relation
\begin{equation}
\label{216}
  [\Lvir_0,J_n({u})]=-n J_n({u})
\end{equation}
for the element $\Lvir_0=J_0(\omega)$, which follows from the axioms for vertex operator algebras. 
Let $\UEA$ be the quotient algebra of the universal enveloping algebra of the Lie algebra $\ggm$ by the two sided ideal generated by $J_0(\vac)-1$ and let us denote the image of $J_n({u})$ by the same symbol. 
We give the degree $d_1+\cdots+d_k$ to the element of the form 
 $
J_{-d_1}({u}_1)\cdots J_{-d_k}({u}_k)
$ 
with ${u}_1,\ldots,{u}_k\in\VOA$. 
Let $\UEA(d)$ be the span of these vectors of degree $d$. 
Then we have
\begin{equation}
  \UEA=\bigoplus_{d\in\ZZ}\UEA(d).
\end{equation}
by which the algebra $\UEA$ becomes a graded algebra.
Note that the relation (\ref{216}) says that the image of $\Lvir_0$ is a Hamiltonian of $\UEA$. 
Consider the standard degreewise topology on $\UEA$ and let $\UEAhat$ denote the degreewise completion. (See Subsection \ref{124}.)
\subsection{The current algebras}
\label{129}
For ${u},{v}\in \VOA$, consider the following expressions in $\UEAhat$ for $k,m,n$:

\begin{equation}
\label{141}
\hskip2em
\eqalign{
\Brm_{k,m,n}({u},{v})
&  =\sum_{i=0}^\infty\binom{k+\Delta({u})-1}{i}J_{k+m+n}({u}_{(n+i)}{v})
\cr&\hskip4em
  -\sum_{i=0}^\infty(-1)^i\binom{n}{i}
   J_{k+n-i}({u})\cdot J_{m+i}({v})
\cr&\hskip6em
  +(-1)^n\sum_{i=0}^\infty(-1)^i\binom{n}{i}
         J_{m+n-i}({v})\cdot J_{k+i}({u}).
\cr&
}
\end{equation}
Then the first sum in the right-hand side is actually a finite sum whereas the second and the last are infinite sums which converge in the linear topology of $\UEAhat(-k-m-n)$. 
The relations $\Brm_{k,m,n}({u},{v})=0$ are nothing else but the Borcherds  identities for the actions on a module where the sums in the right-hand side become finite sums when they act on a fixed element of a module. (See the next subsection for the definition of modules.)
Let $\Bbb$ be the ideal of $\UEAhat$ generated by the elements of the form $\Brm_{k,m,n}({u},{v})$ with ${u},{v}\in \VOA$ and integers $k,m,n$, and let $\hat{\Bbb}$ be the degreewise closure of $\Bbb$. 
Then $\hat{\Bbb}$ is also an ideal of $\UEAhat$.
\begin{remark}
The ideal $\Bbb$ is in fact generated by the elements of the form $\Brm_{k,m,n}({u},{v})$ with $k=-\Delta({u})+1$ and $m,n\in\ZZ$. 
Alternatively, it is also generated by the elements of the form  $\Brm_{k,m,n}({u},{v})$ with $k,m\in\ZZ$ and $n\in\ZZ\setminus\NN$. 
\end{remark}
We now define the {\it current algebra\/} $\UVOA$ associated with  $\VOA$ to be the quotient algebra of $\UEAhat$ by the ideal $\hat{\Bbb}$:
\begin{equation}
  \UVOA=\UEAhat\,/\,\hat\Bbb.
\end{equation}
Then $\UVOA$ is a graded algebra, since $\hat\Bbb$ is a graded ideal, and the image of $\Lvir_0$ is a Hamiltonian. 
\begin{proposition}
The algebra\/ $\UVOA$ is a compatible degreewise complete algebra with Hamiltonian.
\end{proposition}
\begin{note}
This construction is essentially due to Frenkel and Zhu \cite{FZ}. 
The left linear filterwise completion of $\UVOA$ as in Subsection \ref{003} is isomorphic to the current algebra $\mathcal{U}(\VOA)$ considered in \cite{NT}. 
\end{note}
\subsection{Denseness of the current Lie algebra}
Let $\UVOA$ be the current algebra associated with a vertex operator algebra $\VOA$.
Let us regard the current Lie algebra $\ggm$ as a subspace of $\UEA$ and let $\phi$ denote the composition of the canonical maps $\UEA\rightarrow \UEAhat\rightarrow \UVOA$. 
By construction, $\UEA(d)$ is a dense subspace of $\UEAhat(d)$. 
The following observation is insightful.
\begin{proposition}
The image\/ $\phi(\ggm(d))$ is a dense subspace of\/ $\UVOA(d)$ for each integer $d$.
\begin{proof}
It suffices to show that $\phi(\ggm)$ is dense in $\UVOA$ with respect to the left linear topology on $\UVOA$. 
Let us denote by $\phi_n:\UEA\rightarrow \Qrm_n$ the composite of $\phi$ with the canonical surjection $\UVOA\rightarrow \Qrm_n$. 
By the Borcherds relations, we have 
\begin{equation}
\label{229}
\eqalign{
&
J_{s}(u)\cdot J_{t}(v)\cdot 1_n
\cr&
\qquad
=\sum_{m=0}^n\sum_{j=0}^{\Delta(u)+n}(-1)^m\binom{n-s+m}{n-s}\binom{\Delta(u)+n}{j}J_{s+t}(u_{(s+\Delta(u)-m-j-1)}v)\cdot 1_n
\cr
}
\end{equation}
in the quotient $\Qrm_n$ for any integers $s$ and $t$ provided $s\leq n$. 
Hence by induction we have $\phi_n(\UEA)=\phi_n(\ggm)$ for any nonnegative integer $n$. 
Therefore, since $\phi(\UEA(d))$ is dense, $\phi(\ggm(d))$ is also dense. 
\end{proof}
\end{proposition}
\subsection{Exhaustive modules}
Let $\Mmit$ be a vector space and suppose given a series of maps $\VOA\times \Mmit\rightarrow \Mmit$ indexed by integers which we denote by $({u},v)\mapsto \pi^M_n({u})v$. 
Such an $M$ is said to be a {\it weak $\VOA$-module\/} if it satisfies the  conditions listed below.
Set $J^M_n({u})=\pi^M_{n+\Delta({u})-1}({u})$ and let $\Brm^M_{k,m,n}({u},{v})$ be the expression (\ref{141}) with $\{J_n\}$ being replaced by $\{J^M_n\}$.
Then the conditions are as follows:
\vskip2ex
\begin{enumerate}[\rm(i)]
\item
For any ${u}\in\VOA$ and $v\in\Mmit$ there exists an $m$ such that $J^\Mmit_n({u})v=0$ for all $n\geq m$.
%
%
%
\item
The operator $J^\Mmit_n(\vac)$ is the identity if $n=0$ and is zero otherwise.
\item
The identity $\Brm^M_{k,m,n}({u},{v})=0$ holds for any $k,m,n$ and ${u},{v}\in \VOA$.
\end{enumerate}
\vskip2ex
We will say that a weak $\VOA$-module is an {\it exhaustive $\VOA$-module\/} if instead of the condition (i)  the following stronger condition is satisfied:
\vskip2ex
\begin{enumerate}[\rm(i)$'$]
\item
For any $v\in\Mmit$ there exists an $m$ such that $J_{n_1}({u}_1)\cdots J_{n_k}({u}_k)v=0$ for all ${u}_1,\ldots,{u}_k\in \VOA$ whenever $n_1+\cdots+n_k\geq m$. 
\end{enumerate}
\vskip2ex
\begin{remark}
\label{180}
Thanks to the condition (iii), the condition (i)$'$ follows from the apparently weaker condition that for any $v\in\Mmit$ there exists an $m$ such that $J_{n}({u})v=0$ for all ${u}\in \VOA$ and $n\geq m$ by successive use of the relation (\ref{229}).
\end{remark}
Let us consider the map $J_n:\VOA\rightarrow \UVOA$ which sends ${u}\in \VOA$ to the image of $J_n({u})={u}\otimes t^{n+\Delta({u})-1}$ in $\UVOA$.
Then any exhaustive left $\UVOA$-module $M$ becomes an exhaustive $\VOA$-module by letting $J^M_n({u})$ be the action of $J_n({u})$ on $M$.
We will call this $\VOA$-module structure on $M$ the {\it associated $\VOA$-module structure}. 
\begin{proposition}
\label{181}
Let $M$ be an exhaustive $\VOA$-module.
Then there exists a unique structure of an exhaustive\/ $\UVOA$-module on $M$ such that the associated $\VOA$-module structure agrees with the given $\VOA$-module structure on $M$.
\begin{proof}
Let $J^M_n:\VOA\times M\rightarrow M$ be the given $\VOA$-module structure on $M$.
Then they induce a map $\ggm\times M\rightarrow M$ which gives a $\ggm$-module structure on $M$ by the relation $\Brm^M_{k,m,n}({u},{v})=0$ with $n\geq 0$. 
By the universal property of the universal enveloping algebra of $\ggm$, this lifts to a $\UEA$-module structure on $M$ because of the axiom (ii).
Since $M$ is an exhaustive $\VOA$-module, the map $\UEA(d)\times M\rightarrow M$ is continuous for each $d$ when $M$ is endowed with the discrete topology.
Hence this map prolongs to the action of the degreewise completion $\UEAhat$. 
Now the axiom (iii) is nothing else but the defining relations of the algebra $\UVOA$. 
Hence the $\UEAhat$-module structure induces a $\UVOA$-module structure on $M$, which is exhaustive by Lemma \ref{106}.
The uniqueness is clear on each step. 
\end{proof}
\end{proposition}
Thus we have obtained the following result.
\begin{theorem}
\label{142}
The category of exhaustive $\VOA$-modules is canonically equivalent to the category of exhaustive\/ $\UVOA$-modules. 
\end{theorem}
%
%
\section{Associated Poisson algebras}
\label{217}
We will consider the associated graded algebra with respect to a filtration on $\UVOA$ and show that it has a structure of a degreewise complete Poisson algebra. 
\subsection{Zhu's Poisson algebra}
Recall that a (commutative) {\it Poisson algebra\/} is a vector space $\Poisson$  equipped with two bilinear maps $\cdot:\Poisson\times \Poisson\rightarrow \Poisson$ and $\{\ ,\ \}:\Poisson\times \Poisson\rightarrow \Poisson$ called the multiplication and the Poisson bracket, respectively, such that $\Poisson$ is a commutative associative algebra with unity with respect to the multiplication, $\Poisson$ is a Lie algebra with respect to the Poisson bracket and the Leibniz identity holds:
\begin{equation}
\{{x}\cdot{y},{z}\}={x}\cdot\{{y},{z}\}+{y}\cdot\{{x},{z}\}
.
\end{equation}
We denote the unity of $\Poisson$ by $1_\Poisson$. 
Let $\VOA$ be a vertex operator algebra. 
We let $\Crm_2(\VOA)$ be the span of the elements of the form ${u}_{(n)}{v}$ with ${u},{v}\in \VOA$ and $n\leq -2$. 
We set
\begin{equation}
  {u}\cdot {v}={u}_{(-1)}{v}
\quad\hbox{and}\quad
  \{{u},{v}\}={u}_{(0)}{v}.
\end{equation}
The following result is obtained in \cite{Zhu}.
\begin{proposition}[Zhu]
The operations $\cdot$ and $\{\ ,\ \}$ induces a Poisson algebra structure on the quotient space $\VOA/\Crm_2(\VOA)$.
\end{proposition}
Let us denote by $\Poisson$ the Poisson algebra $\VOA/\Crm_2(\VOA)$ obtained as above, which we will call {\it Zhu's Poisson algebra}.

\subsection{Poisson filtrations and the associated graded algebras}
\label{130}
Let $\VOA$ be a vertex operator algebra and let $\ggm$, $\UEA$, $\UEAhat$ and $\UVOA$ as in the preceding section.
Let $\Grm_p\ggm$ be the image of $\bigoplus_{r\leq p}\VOA[r]\otimes_\KK\KK[t,t^{-1}]$ in $\ggm$ and let $\Grm_p\UEA$ be the sum of subspaces $\Grm_{p_1}\ggm\cdots\Grm_{p_k}\ggm$ with $k=0,1,2\ldots$ and $p_1+\cdots+p_k=p$ in $\UEA$. 
Then $\Grm$ is a separated filtration on $\UEA$ satisfying
\begin{equation}
\label{173}
  \Grm_p\UEA\cdot \Grm_q\UEA\subset \Grm_{p+q}\UEA
\end{equation}
for any integers $p$ and $q$.
Let $\Grm_p\UEAhat(d)$ be the closure of the image of $\Grm_p\UEA(d)=\Grm_p\UEA\cap \UEA(d)$ in $\UEAhat(d)$. 
Then the associated graded algebra is given by 
\begin{equation}
  \gr^\Grm\UEAhat=\bigoplus_d\bigoplus_p\gr^\Grm_p\UEAhat(d)
,\quad 
  \gr^\Grm_p\UEAhat(d)=\Grm_p\UEAhat(d)\,/\,\Grm_{p-1}\UEAhat(d)
.
\end{equation}
Considering the quotients by $\hat\Bbb$, we obtain the induced filtration $\Grm_p\UVOA$ of the current algebra $\UVOA$ and the associated graded algebra: 
\begin{equation}
  \gr^\Grm\UVOA
    =\bigoplus_{d=-\infty}^\infty \gr^\Grm\UVOA(d)
,\quad
  \gr^\Grm\UVOA(d)
    =\bigoplus_{p} \gr^\Grm_p\UVOA(d)
.
\end{equation}
Let us give the space $\gr^\Grm\UVOA$ the induced degreewise topology and let $\SPoisson$ be the degreewise completion of the algebra $\gr^\Grm\UVOA$:
\begin{equation}
  \SPoisson=\bigoplus_{d}\SPoisson(d)
,\quad
  \SPoisson(d)=\varprojlim_{n}\gr^\Grm\UVOA(d)\!\bigm/\!\Igraded_n(\gr^\Grm\UVOA(d))
.
\end{equation}
Then this is a compatible degreewise complete algebra. 
We will denote the image of $J_n({u})$ in $\SPoisson$ by $\Psi_n({u})$. 
\begin{proposition}
\label{187}
The algebra\/ $\UVOA$ is quasi-finite if and only if\/ $\SPoisson$ is so.
\begin{proof}
By the construction, we have 
\begin{equation}
  \Qrm_n(\SPoisson(d))=\gr^\Grm\Qrm_n(\UVOA(d))
,\quad
  \Qrm_n(\UVOA(d))=\bigcup_p \Grm_p\Qrm_n(\UVOA(d))
,
\end{equation}
where $\Grm_p\Qrm_n(\UVOA(d))$ denotes the induced filtration. 
Hence $\Qrm_n(\UVOA(d))$ is finite-dimensional if and only if $\Qrm_n(\SPoisson(d))$ is so. 
\end{proof}
\end{proposition}
\subsection{Associated Poisson structure}
Consider the operation of taking commutator of elements of $\UVOA$:
\begin{equation}
  \UVOA\times \UVOA\rightarrow \UVOA
,\quad
  (a,b)\mapsto [a,b]=a\cdot b-b\cdot a.
\end{equation}
The particular case of the identities $\Brm_{k,m,n}({u},{v})=0$ with $k=0$ is given by 
\begin{equation}
\label{189}
  [J_m({u}),J_n({v})]
    =\sum_{i=0}^\infty \binom{m+\Delta({u})-1}{i} J_{m+n}({u}_{(i)}{v})
.
\end{equation}
\begin{lemma}
\label{172}
 $
  [\Grm_p\UVOA,\Grm_q\UVOA]\subset \Grm_{p+q-1}\UVOA
$. 
\begin{proof}
The left-hand side of (\ref{189}) belongs to
 $\Grm_{\Delta({u})+\Delta({v})}\UVOA$ whereas the element
 $J_{m+n}({u}_{(i)}{v})$ in the right-hand side belongs to
 $\Grm_{\Delta({u})+\Delta({v})-i-1}\UVOA$ for $i\geq 0$. 
Hence an element of $[\Grm_p\UVOA,\Grm_q\UVOA]$ is written as a sum of elements of $\Grm_{p+q-1}\UVOA$. 
\end{proof}
\end{lemma}
Thanks to this lemma, the multiplication of $\gr^\Grm\UVOA$ is commutative and the operation $[\ ,\ ]:\Grm_p\UVOA\times \Grm_q\UVOA\rightarrow \Grm_{p+q-1}\UVOA$ induces an operation
\begin{equation}
  \gr^\Grm_p\UVOA\times \gr^\Grm_q\UVOA\rightarrow \gr^\Grm_{p+q-1}\UVOA
,\quad
  (\alpha,\beta)\mapsto \{\alpha,\beta\}
\end{equation}
by letting $\{\alpha,\beta\}$ be the image of $[a,b]$ in $\gr^\Grm_{p+q-1}\UVOA$, where $a$ and $b$ are representatives of $\alpha$ and $\beta$, respectively. 

In general, we will call a compatible degreewise topological algebra with a Poisson algebra structure for which the Poisson bracket is continuous a {\it compatible degreewise topological Poisson algebra}. 
In case the degreewise topology is complete then we will say that the Poisson algebra is a {\it compatible degreewise complete Poisson algebra}. 
\begin{proposition}
\label{174}
The multiplication and the bracket operation defined as above endow the space\/ $\SPoisson$ with a structure of a compatible degreewise complete Poisson algebra.
%
%
\end{proposition}
\subsection{Relation to Zhu's Poisson algebra}
Let us look more carefully at the relations $\Brm_{k,m,n}({u},{v})=0$. 
When $k=-\Delta({u})+1$ and $m+n=p+\Delta({u})-1$, we have 
\begin{equation}
\label{178}
\eqalign{
&  J_{p}({u}_{(n)}{v})
%
%
  =\sum_{i=0}^\infty(-1)^i\binom{n}{i}J_{k+n-i}({u})\cdot J_{m+i}({v})
\cr&\hskip6em
  -(-1)^n\sum_{i=0}^\infty(-1)^i\binom{n}{i}
         J_{m+n-i}({v})\cdot J_{k+i}({u}).
\cr&
}
\end{equation}
Then the left-hand side belongs to $\Grm_{\Delta({u})+\Delta({v})-n-1}\UVOA$ whereas the right-hand side to $\Grm_{\Delta({u})+\Delta({v})}\UVOA$. 
Therefore, 
\begin{equation}
\label{188}
  J_p({u}_{(n)}{v})\equiv 0\hbox{\quad if $n\leq -2$}
.
\end{equation}
This implies that the map $J_p:\VOA\rightarrow \gr^\Grm\UVOA$ factors a map 
\begin{equation}
\Psi_p:\Poisson\rightarrow\gr^\Grm\UVOA
\end{equation}
where $\Poisson$ is Zhu's Poisson algebra $\VOA/\Crm_2(\VOA)$.

Now substitute $n=-1$ in (\ref{178}).
Then we have
\begin{equation}
\label{214}
\eqalign{
&  J_{p}({u}_{(-1)}{v})
  =\sum_{i=0}^\infty J_{k-i-1}({u})\cdot J_{m+i}({v})
  +\sum_{i=0}^\infty J_{m-i-1}({v})\cdot J_{k+i}({u}).
\cr&}
\end{equation}
Projecting this relation to the associated graded algebra, we immediately have
\begin{equation}
\label{215}
\eqalign{
&  \Psi_{p}({u}_{(-1)}{v})
  =\sum_{j=-\infty}^\infty \Psi_{p-j}({u})\cdot \Psi_{j}({v}).
\cr&
}
\end{equation}

\begin{note}
The results in this section and the next are reformulations of the arguments in Subsection 3.2 of \cite{NT}. 
\end{note}
%
\section{Poisson current algebras}
\label{027}
In this section, we will construct a universal Poisson algebra satisfying the relation (\ref{215}), which we will call a Poisson current algebra, and will investigate its properties. 
\subsection{Symmetric algebras on the loop Lie algebras}
\label{167}
Let $\Poisson$ be a Poisson algebra.
A {\it Poisson ideal\/} of $\Poisson$ means a subspace $\agm$ of $\Poisson$ such that both $\Poisson\cdot \agm\subset \agm$ and $\{\Poisson,\agm\}\subset \agm$ hold.
Consider the case when $\Poisson$ is given a grading $\Poisson=\bigoplus_{r}\Poisson[r]$ indexed by integers satisfying
\begin{equation}
\label{165}
  \{\Poisson[r],\Poisson[s]\}\subset \Poisson[r+s-1]
,\quad 
\Poisson[r]\cdot \Poisson[s]\subset \Poisson[r+s].
\end{equation}
Then the unity $1_\Poisson$ must belong to $\Poisson[0]$.
We will call a Poisson algebra endowed with such a grading a {\it graded Poisson algebra}.
We denote $\Delta({x})=r$ when ${x}\in \Poisson[r]$ and call it the weight of ${x}$.
Recall the well-known fact that the symmetric algebra on a Lie algebra has a canonical structure of a Poisson algebra induced from the Lie bracket operation on the Lie algebra. 
Let $\Poisson=\bigoplus_r\Poisson[r]$ be a graded Poisson algebra and let $\Poisson[t,t^{-1}]$ be the loop Lie algebra
\begin{equation}
  \Poisson[t,t^{-1}]=\Poisson\otimes_{\KK}\KK[t,t^{-1}]
\end{equation}
with the Lie bracket defined by
\begin{equation}
  [{x}\otimes t^m,{y}\otimes t^n]=\{{x},{y}\}\otimes t^{m+n}.
\end{equation}
For each homogeneous ${x}$, we set 
\begin{equation}
  \tilde{\Psi}_n({x})={x}\otimes t^{n+\Delta({x})-1}
\end{equation}
and extend it linearly to all ${x}$. 
Then the Lie bracket operation takes the following form: 
\begin{equation}
  [\tilde{\Psi}_m({x}),\tilde{\Psi}_n({y})]=\tilde{\Psi}_{m+n}(\{{x},{y}\}).
\end{equation}
Let $\qgm$ be the quotient of $\Poisson[t,t^{-1}]$ by the span of the elements  $\tilde{\Psi}_n(1_\Poisson)$ with $n\ne {0}$. 
We will denote the image of $\tilde{\Psi}_n({x})$ in $\qgm$ by the same symbol.
Since $[\tilde{\Psi}_m(1_\Poisson),\tilde{\Psi}_n({x})]=0$, the space $\qgm$ becomes a Lie algebra. 
Now let $\Sym$ be the quotient algebra of the symmetric algebra on $\qgm$ by the ideal generated by $\tilde{\Psi}_0(1_\Poisson)-1$ and let us denote the image of $\tilde{\Psi}_n(x)$ in $\Sym$ by the same symbol. 
We give the degree $d_1+\cdots+d_k$ to the element of the form 
 $
\tilde{\Psi}_{-d_1}({x}_1)\cdots\tilde{\Psi}_{-d_k}({x}_k)
$ 
with ${x}_1,\ldots,{x}_k\in\Poisson$. 
Let $\Sym(d)$ be the span of these vectors of degree $d$. 
Then we have 
\begin{equation}
  \Sym=\bigoplus_d\Sym(d)
\end{equation}
and the algebra $\Sym$ becomes a Poisson algebra. Note that we have
\begin{equation}
  \Sym(d)\cdot \Sym(e)\subset \Sym(d+e)
,\quad
  \{\Sym(d),\Sym(e)\}\subset \Sym(d+e)
.
\end{equation}
Recall the standard degreewise topology on $\Sym$ defined by
\begin{equation}
  \Igraded_n(\Sym(d))=\sum_{k\leq -n-1}\Sym(d-k)\cdot \Sym(k), 
\end{equation}
which is separated. 
Then the multiplication maps $\Sym(d)\times\Sym(e)\rightarrow \Sym(d+e)$ is  continuous. 
Moreover we have the following.
\begin{lemma}
\label{159}
The Poisson bracket operation\/ $\Sym(d)\times \Sym(e)\rightarrow \Sym(d+e)$ is continuous with respect to the standard degreewise topology.
\begin{proof}
Let $i$ be any integer with $i\leq -n-1$.
Then we have
$$
\hskip2em
\eqalign{
&\hskip-2em
  \{\Sym(d),\Sym(e-i)\cdot \Sym(i)\}
\cr&
    \subset \{\Sym(d),\Sym(e-i)\}\cdot \Sym(i)
             +\{\Sym(d),\Sym(i)\}\cdot \Sym(e-i)
\cr&
    =        \Sym(d+e-i)\cdot \Sym(i)
             +\Sym(e-i)\cdot \Sym(d+i)
.
\cr
}
$$
Therefore, $\{{u}+\Igraded_k(\Sym(d)),{v}+\Igraded_m(\Sym(e))\}\subset \{{u},{v}\}+\Igraded_n(\Sym(d+e))$ if $k$ and $m$ satisfy $k,k-e,m,m-d\geq n$. 
\end{proof}
\end{lemma}
Let $\hat{\Sym}$ be the degreewise completion of $\Sym$ with respect to the standard degreewise topology. 
Then by Lemma \ref{159} the Poisson algebra structure on $\Sym$ extends to $\hat{\Sym}$.
Let us denote the image of $\tilde{\Psi}_n({x})$ under the canonical map $\Sym\rightarrow \hat{\Sym}$ again by the same symbol.
The algebra $\hat\Sym$ is a compatible degreewise complete Poisson algebra. 
\subsection{Poisson current algebras}
\label{030}
For each ${x},{y}\in\Poisson$ and an integer $m$, we set
\begin{equation}
\eqalign{
&
  \Drm_m({x},{y})
  =\tilde{\Psi}_m({x}\cdot {y})
      -\sum_{j\in\ZZ}\tilde{\Psi}_{m-j}({x})\cdot \tilde{\Psi}_j({y})
.\cr}
\end{equation}
and let $\hat{\Dbb}$ denote the degreewise closure of the ideal of $\hat{\Sym}$ generated by the elements $\Drm_m({x},{y})$ with ${x},{y}\in\Poisson$ in $\hat{\Sym}$ and $m\in\ZZ$.
\begin{lemma}
\label{158}
The ideal\/ $\hat{\Dbb}$ is a Poisson ideal of\/ $\hat{\Sym}$.
\begin{proof}
Recall that the Poisson bracket operation is continuous.
Hence we have
$$
\hskip2em
\eqalign{
&\hskip-2em
  \Bigl\{\Drm_m({x},{y}),\tilde{\Psi}_n({z})\Bigr\}
\cr&
    =\tilde{\Psi}_{m+n}(\{{x}\cdot{y},{z}\})
    -\sum_j
      \{\tilde{\Psi}_{m-j}({x})\cdot \tilde{\Psi}_j({y}),\tilde{\Psi}_n({z})\}
\cr&
    =\tilde{\Psi}_{m+n}({x}\cdot \{{y},{z}\})
    -\sum_j
      \tilde{\Psi}_{m-j}({x})\cdot \tilde{\Psi}_{j+n}(\{{y},{z}\})
\cr&\qquad 
    +\tilde{\Psi}_{m+n}({y}\cdot \{{x},{z}\})
    -\sum_j
      \tilde{\Psi}_j({y})\cdot \tilde{\Psi}_{m+n-j}(\{{x},{z}\})
\cr&
    =\Drm_{m+n}({x},\{{y},{z}\})+\Drm_{m+n}({y},\{{x},{z}\})
,
\cr
}$$
so the conclusion follows. 
\end{proof}
\end{lemma}
We let $\SPoissontilde=\SPoissontilde(\Poisson)$ be the quotient of $\hat\Sym$ by the Poisson ideal $\hat{\Dbb}$:
\begin{equation}
  \SPoissontilde=\hat\Sym\,/\,\hat\Dbb
.
\end{equation}
We will call the Poisson algebra $\SPoissontilde$ the {\it Poisson current algebra\/} associated with the Poisson algebra $\Poisson$.
\begin{proposition}
The algebra\/ $\SPoissontilde$ is a compatible degreewise complete Poisson algebra.
\end{proposition}
\subsection{Quasi-finiteness}
\label{058}
Now assume that our Poisson algebra $\Poisson$ is finite-dimensional and let ${x}_1,\ldots,{x}_r$ be a basis of a linear complement of $\KK 1_\Poisson$ in $\Poisson$. 
Then the algebra $\Sym$ is spanned by the elements of the form $\tilde\Psi_{-d_1}({x}_{i_1})\cdots\tilde\Psi_{-d_k}({x}_{i_k})$ with $k\geq 0$ and $d_1\geq \cdots \geq d_k$. 
Consider the left canonical quotient module $\Qrm_n$ in the case when the algebra $\Agraded$ is $\SPoissontilde$.  
The module $\Qrm_n$ is generated as a left $\SPoissontilde$-module by the image $1_{n}$ of the unit of $\SPoissontilde$. 
For each $k\geq 0$ and $d\geq -n$, consider the set
\begin{equation}
  \Pi_k(d)=\ensemble{(d_1,\ldots,d_k)}%
               {$d_1\geq \cdots\geq d_k\geq -n$ and $d_1+\cdots+d_k=d$}.
\end{equation}
We will call a vector of the form $\tilde\Psi_{-d_1}({x}_{i_1})\cdots\tilde\Psi_{-d_k}({x}_{i_k})\cdot 1_{n}$ a {\it vector with index $(d_1,\ldots,d_k)$}. 
We put
\begin{equation}
  \interior{\Pi}_k(d)=\ensemble{(d_1,\ldots,d_k)}%
               {$d_1>\cdots>d_k\geq -n$ and $d_1+\cdots+d_k=d$}
\end{equation}
and set 
\begin{equation}
  {\Pi}(d)=\bigcup_{k=0}^\infty {\Pi}_k(d)
,\quad
  \interior{\Pi}(d)=\bigcup_{k=0}^\infty \interior{\Pi}_k(d)
.
\end{equation}
\begin{lemma}
\label{161}
The space\/ $\Qrm_n(d)$ is spanned by vectors with indices in\/ $\interior{\Pi}(d)$.
\begin{proof}
Since the module $\Qrm_n$ is exhaustive and the image of each $\Sym(d)$ is dense in $\SPoissontilde(d)$, the space $\Qrm_n(d)$ is spanned by vectors with indices in $\Pi(d)$. 
Let us show that any vector with index in $\Pi_k(d)$ is a linear combination of vectors with indices in $\interior{\Pi}(d)$. 
Introduce the lexicographic order on the set $\Pi_k(d)$: we define $(d_1,\ldots,d_k)<(e_1,\ldots,e_k)$ by $d_1<e_1$ and in case $d_1=e_1$  recursively by $(d_2,\ldots,d_k)<(e_2,\ldots,e_k)$. 
Then as this is a total order on a finite set, there exists a maximum element: that is $(d+(k-1)n,-n,\ldots,-n)$. 
We now proceed by induction on the length $k$. 
The case $k=1$ is trivial. 
Assume that the claim is true for any vector with index shorter than $k$ and suppose given a vector
\begin{equation}
\label{160}
  \tilde\Psi_{-d_1}({x}_{i_1})\cdots\tilde\Psi_{-d_k}({x}_{i_k})\cdot 1_{n}
\end{equation}
with $d_1\geq\cdots\geq d_k>-n-1$.
If $d_1>\cdots>d_k>-n-1$ then we have nothing to prove so we consider the case when $d_i=d_{i+1}$ at some position $i$.
Recall the relation
 $
  \tilde\Psi_{2k}({x}\cdot {y})=\sum_j \tilde\Psi_{k-j}({x})\cdot \tilde\Psi_{k+j}({y})
$ 
which implies
\begin{equation}
\label{171}
  \tilde\Psi_{-d_i}({x}_i)\cdot \tilde\Psi_{-d_i}({x}_{i+1})=\tilde\Psi_{-2d_i}({x}_i\cdot {x}_{i+1})-\sum_{j\ne 0}\tilde\Psi_{-d_i-j}({x}_i)\cdot \tilde\Psi_{-d_i+j}({x}_{i+1})
.
\end{equation}
Hence the vector (\ref{160}) is rewritten as the sum of a shorter vector and a finite number of vectors greater in the lexicographic order. 
By the inductive hypothesis, the shorter vector is written by the vectors with indices in $\interior{\Pi}(d)$. 
Then apply the same argument to the vectors with greater indices in the rest of the sum. 
The recursion stops within a finite number of steps, at most at the maximum.
\end{proof}
\end{lemma}
\begin{note}
The argument described above is a refined variation of the proof of Theorem 3.2.7 in \cite{NT}. 
The idea of utilizing the relation (\ref{171}) goes back to \cite{GN}. 
See \cite{Bu} and \cite{Li2} for related results. 
\end{note}
We will say that a compatible degreewise complete Poisson algebra is {\it quasi-finite\/} if the conditions of quasi-finiteness for a compatible degreewise complete algebras are satisfied except the existence of a Hamiltonian. 
Now the following result is an immediate consequence of Lemma \ref{161}.
\begin{theorem}
\label{177}
If the Poisson algebra\/ $\Poisson$ is finite-dimensional then the algebra\/ $\SPoissontilde$ is quasi-finite.
\begin{proof}
By Lemma \ref{161}, the space $\Qrm_n(d)$ is spanned by the vectors associated with $\interior{\Pi}(d)=\bigcup_{k=0}^\infty\interior{\Pi}_k(d)$, which is a finite set. 
\end{proof}
\end{theorem}
%
\section{Current algebras and Poisson current algebras}
%
In this section, we will show that Zhu's finiteness condition on a vertex operator algebra implies the quasi-finiteness of the associated current algebra.
This will be done by relating the results of the preceding section to the  Poisson algebra $\SPoisson$ associated with $\UVOA$. 
\subsection{Relation to the current algebras}
A {\it homomorphism of degreewise topological Poisson algebras\/} is a map from a degreewise topological Poisson algebra to another such that it is a homomorphism of Poisson algebras that preserves the gradings for which the restriction to each homogeneous subspace is continuous. 

Let $\VOA$ be a vertex operator algebra and let $\Poisson$ be Zhu's Poisson algebra $\VOA/\Crm_2(\VOA)$.
Recall the notations in the previous sections of Part II. 

By the definition of Zhu's Poisson algebra, (\ref{188}) implies that the maps $\VOA\rightarrow \gr^\Grm\UVOA$ which sends ${u}\in \VOA[r]$ to the image $\Psi_p({u})$ of $J_p({u})$ in $\gr^\Grm_r\UVOA$ factors a map
\begin{equation}
\label{190}
  \Psi_p:\Poisson\rightarrow \gr^\Grm\UVOA
.
\end{equation}
Then the set of the maps $\Psi_p:\Poisson\rightarrow \gr^\Grm\UVOA$ gives rise to a single map 
\begin{equation}
  \Psi:\qgm\rightarrow \gr^\Grm\UVOA
\end{equation}
which sends $\tilde{\Psi}_n({x})$ to $\Psi_n({u})$, where ${x}=\bar{{u}}$ is the class of ${u}\in \VOA$ in $\Poisson$.
Now the relations (\ref{189}) and (\ref{188}) imply that the map $\Psi$ is a homomorphism of Lie algebras. 
Since $\gr^\Grm\UVOA$ is a Poisson algebra, this map induces a unique homomorphism
\begin{equation}
  \Sym\rightarrow \gr^\Grm\UVOA
\end{equation}
of Poisson algebras by the universal property of the symmetric algebra.
We denote this map by the same symbol $\Psi$.
\begin{lemma}
\label{203}
The map\/ $\Psi$ prolongs to a surjective homomorphism\/ $\hat{\Sym}\rightarrow \SPoisson$ of degreewise topological Poisson algebras. 
\begin{proof}
The assertion follows immediately from the construction by noting the relation  $\Psi(\Igraded_n(\Sym))=\Igraded_n(\gr^\Grm\UVOA)$, where the latter space is the one induced from $\Igraded_n(\UVOA)$. 
\end{proof}
\end{lemma}
Now let $\SPoissontilde$ be the Poisson current algebra of $\Poisson$ as defined in Subsection \ref{030}. 
The relation (\ref{215}) implies that the ideal $\hat{\Dbb}$ is mapped to the closure of $\gr^\Grm\hat{\Bbb}$.
Therefore, the map $\hat{\Sym}\rightarrow \SPoisson$ induces a homomorphism $\SPoissontilde\rightarrow \SPoisson$ of degreewise topological Poisson algebras. 
Thus we have verified the following result.
\begin{theorem}
\label{176}
Let $\VOA$ be a vertex operator algebra and let\/ $\UVOA$ be the associated current algebra.
Let $\SPoisson$ be the degreewise completion of $\gr^\Grm\UVOA$ and let\/ $\SPoissontilde$ be the Poisson current algebra associated with Zhu's Poisson algebra $\Poisson=\VOA/\Crm_2(\VOA)$. 
Then there exists a surjective homomorphism\/ $\SPoissontilde\rightarrow \SPoisson$ of degreewise topological Poisson algebras.
\end{theorem}
\subsection{Consequences of Zhu's finiteness condition}
A vertex operator algebra $\VOA$ is said to satisfy {\it Zhu's finiteness condition\/} or said to be {\it $\mathit{C}_{\mathit 2}$-finite\/} if Zhu's Poisson algebra $\Poisson=\VOA/\Crm_2(\VOA)$ is finite-dimensional. 

By combining Proposition \ref{187}, Theorem \ref{176} and Theorem \ref{177}, we immediately see that Zhu's finiteness condition implies quasi-finiteness. 
Namely, we have the following theorem which is the main result of Part II. 
\begin{theorem}
\label{136}
If a vertex operator algebra $\VOA$ satisfies Zhu's finiteness condition then the associated current algebra\/ $\UVOA$ is quasi-finite. 
\end{theorem}

Let $\VOA$ be a vertex operator algebra satisfying Zhu's finiteness condition and let $\UVOA$ be the associated current algebra. 
Let $\UVOA_n$ be the finite-dimensional algebra associated with $\UVOA$ as 
defined in Subsection \ref{100} and let $\gap$ be the number defined by (\ref{231}).
Then Theorem \ref{136} allows us to apply the results of Part I to the algebra $\Agraded=\UVOA$. 
For instance, we have the following.

\begin{corollary}
Let $\VOA$ be a $C_2$-finite vertex operator algebra and let $n$ be an integer such that $n\geq \gap$. 
Then the category of exhaustive $\VOA$-modules is canonically equivalent to the category of left\/ $\UVOA_n$-modules. 
\end{corollary}

Moreover, the finiteness theorems in Subsection \ref{224} hold for the various categories of $\UVOA$-modules. 
%
%
%
\def\refname
\end{document}